\documentclass{article} 
\usepackage{iclr2026_conference,times}

\usepackage{hyperref}
\definecolor{mydarkblue}{rgb}{0,0.08,0.45}
\hypersetup{ %
  colorlinks=true,
  linkcolor=mydarkblue,
  citecolor=mydarkblue,
  filecolor=mydarkblue,
  urlcolor=mydarkblue,
  }
\usepackage{url}

\title{
    {\red Ringleader ASGD}: The First Asynchronous SGD
    with Optimal Time Complexity under Data Heterogeneity
    }

\author{Artavazd Maranjyan, Peter Richt\'{a}rik \\
King Abdullah University of Science and Technology (KAUST) \\
\texttt{\{arto.maranjyan,richtarik\}@gmail.com} \\
}

\usepackage{graphicx}
\usepackage{apptools}
\usepackage[flushleft]{threeparttable}
\usepackage{array,booktabs,makecell}
\usepackage{multirow}
\usepackage{nicefrac}
\usepackage{amsthm}
\usepackage{mathtools}
\usepackage{amsmath,amsfonts,bm,amssymb}
\usepackage{dsfont}
\usepackage{wrapfig}
\usepackage{caption}
\usepackage{siunitx}
\usepackage{thm-restate}
\usepackage{nccmath}
\usepackage{empheq}
\usepackage{bbm}
\usepackage{tabularx}

\usepackage{algorithmic}
\usepackage{algorithm}
\usepackage{filecontents}

\usepackage[capitalize,noabbrev]{cleveref}

\usepackage{color}
\usepackage{colortbl}
\definecolor{bgcolor}{rgb}{0.76,0.88,0.50}
\definecolor{bgcolor0}{rgb}{0.93,0.99,1}
\definecolor{bgcolor1}{rgb}{0.8,1,1}
\definecolor{bgcolor2}{rgb}{0.8,1,0.8}
\definecolor{bgcolor3}{rgb}{0.50,0.90,0.50}
\usepackage{tcolorbox}
\usepackage{pifont}

\definecolor{mydarkgreen}{RGB}{39,130,67}
\definecolor{mydarkorange}{RGB}{236,147,14}
\definecolor{mydarkred}{RGB}{192,47,25}
\definecolor{ruby}{RGB}{155,17,30}
\definecolor{chili}{RGB}{191,0,0}
\definecolor{sangria}{RGB}{146,0,10}
\definecolor{burgundy}{RGB}{128,0,32} 
\definecolor{darkred}{RGB}{132,0,0} 
\definecolor{cherry}{RGB}{192,0,0} 
\definecolor{blue}{RGB}{0,0,255}

\newcommand{\orange}{\color{mydarkorange}}

\newcommand{\red}{\color{cherry}}

\usepackage[textsize=tiny]{todonotes}

\usepackage{xspace}
\usepackage{relsize}
\newcommand{\algname}[1]{{\sf\relscale{0.95}#1}\xspace}

\newcommand{\norm}[1]{\left\| #1 \right\|}
\newcommand{\sqnorm}[1]{\left\| #1 \right\|^2}

\newcommand{\inp}[2]{\left\langle#1,#2\right\rangle} 

\newcommand{\R}{\mathbb{R}} 
\newcommand{\E}[1]{\mathbb{E}\left[#1\right]}

\newcommand{\Exp}[1]{{\mathbb{E}}\left[#1\right]}
\newcommand{\ExpSub}[2]{{\mathbb{E}}_{#1}\left[#2\right]}
\newcommand{\ExpCond}[2]{{\mathbb{E}}\left[\left.#1\ \right\vert\ #2\right]}

\newcommand{\cD}{\mathcal{D}}

\newcommand{\cN}{\mathcal{N}}
\newcommand{\cO}{\mathcal{O}}

\newcommand{\eqdef}{:=} 

\newcommand{\minimize}{\mathop{\mathrm{minimize}}}

\newcommand{\tauavg}{\tau_{\mathrm{avg}}}

\makeatletter
\newcommand{\vast}{\bBigg@{4}}

\def\<{\left\langle}
\def\>{\right\rangle}
\def\[{\left[}
\def\]{\right]}
\def\({\left(}
\def\){\right)}

\definecolor{marks}{rgb}{0.3,0.25,0.2}
\newcommand{\checkmarkgreen}{\textbf{\color{marks}\ding{52}}}
\newcommand{\crossmarkred}{{\color{marks}\ding{56}}}

\usepackage{thmtools}
\usepackage{thm-restate}
\usepackage{mdframed}

\theoremstyle{plain}
\newtheorem{theorem}{Theorem}[section]

\newtheorem{remark}[theorem]{Remark}

\theoremstyle{definition} 

\newmdenv[
  font=\normalfont\bfseries,
  linecolor=black,
  linewidth=0.8pt,
  topline=false,
  bottomline=false,
  leftline=false,
  rightline=false,
  backgroundcolor=gray!13,
  skipabove=10pt,
  skipbelow=10pt,
  innertopmargin=5pt,
  innerbottommargin=4pt,
  innerrightmargin=4pt,
  innerleftmargin=4pt,
]{myboxed}

\newmdtheoremenv[
  font=\normalfont\bfseries,
  linecolor=black,
  linewidth=0.8pt,
  topline=false,
  bottomline=false,
  leftline=false,
  rightline=false,
  backgroundcolor=gray!13,
  skipabove=10pt,
  skipbelow=10pt,
  innertopmargin=5pt,
  innerbottommargin=4pt,
  innerrightmargin=4pt,
  innerleftmargin=4pt,
]{boxedtheorem}{Theorem}

\newmdtheoremenv[
  font=\normalfont\bfseries,
  linecolor=black,
  linewidth=0.8pt,
  topline=false,
  bottomline=false,
  leftline=false,
  rightline=false,
  backgroundcolor=gray!13,
  skipabove=4pt,
  skipbelow=12pt,
  innertopmargin=5pt,
  innerbottommargin=4pt,
  innerrightmargin=4pt,
  innerleftmargin=4pt,
]{boxedlemma}{Lemma}

\newmdtheoremenv[
  font=\normalfont\bfseries,
  linecolor=black,
  linewidth=0.8pt,
  topline=false,
  bottomline=false,
  leftline=false,
  rightline=false,
  backgroundcolor=gray!13,
  skipabove=10pt,
  skipbelow=10pt,
  innertopmargin=5pt,
  innerbottommargin=4pt,
  innerrightmargin=4pt,
  innerleftmargin=4pt,
]{boxeddefinition}{Definition}

\newmdtheoremenv[
  font=\normalfont\bfseries,
  linecolor=black,
  linewidth=0.8pt,
  topline=false,
  bottomline=false,
  leftline=false,
  rightline=false,
  backgroundcolor=gray!13,
  skipabove=10pt,
  skipbelow=10pt,
  innertopmargin=5pt,
  innerbottommargin=4pt,
  innerrightmargin=4pt,
  innerleftmargin=4pt,
]{boxedassumption}{Assumption}

\newtheorem{innercustomthm}{Theorem}
\newenvironment{restate-theorem}[1]
  {\renewcommand\theinnercustomthm{#1}\innercustomthm}
  {\endinnercustomthm}

\newtheorem{innercustomlemma}{Lemma}
\newenvironment{restate-lemma}[1]
  {\renewcommand\theinnercustomlemma{#1}\innercustomlemma}
  {\endinnercustomlemma}

\newtheorem{innercustomproposition}{Proposition}
\newenvironment{restate-proposition}[1]
  {\renewcommand\theinnercustomproposition{#1}\innercustomproposition}
  {\endinnercustomproposition}

\newenvironment{restate-boxedtheorem}[1]
  {\renewcommand\theinnercustomthm{#1}\begin{myboxed}\begin{innercustomthm}}
  {\end{innercustomthm}\end{myboxed}}

\newenvironment{restate-boxedlemma}[1]
  {\renewcommand\theinnercustomlemma{#1}\begin{myboxed}\begin{innercustomlemma}}
  {\end{innercustomlemma}\end{myboxed}}

\newcommand*{\sketchproofname}{Sketch of Proof}

\newcommand{\algn}{{\sf\red\relscale{0.95}Ringleader ASGD}\xspace}
\newcommand{\iasgd}{{\sf\relscale{0.95}IA\textsuperscript{2}SGD}\xspace}
\newcommand{\iasgdtitle}{IA\textsuperscript{2}SGD}
\newcommand{\malenia}{\algname{Malenia SGD}}
\newcommand{\maleniatitle}{Malenia SGD}
\newcommand{\minibatch}{\algname{Minibatch SGD}}

\newcommand{\naiveminibatch}{\algname{Naive Minibatch SGD}}
\newcommand{\naiveminibatchtitle}{Naive Minibatch SGD}

\newcommand{\saga}{\algname{SAGA}}

\newcommand{\ringmaster}{\algname{Ringmaster ASGD}}

\newcommand{\rennala}{\algname{Rennala SGD}}
\newcommand{\sgd}{\algname{SGD}}

\iclrfinalcopy 
\begin{document}
\maketitle
\begin{abstract}

    Asynchronous stochastic gradient methods are central to scalable distributed optimization, particularly when devices differ in computational capabilities.
    Such settings arise naturally in federated learning, where training takes place on smartphones and other heterogeneous edge devices.
    In addition to varying computation speeds, these devices often hold data from different distributions.
    However, existing asynchronous \sgd methods struggle in such heterogeneous settings and face two key limitations.
    First, many rely on unrealistic assumptions of similarity across workers' data distributions.
    Second, methods that relax this assumption still fail to achieve theoretically optimal performance under heterogeneous computation times.
    We introduce \algn, the first asynchronous \sgd algorithm that attains the theoretical lower bounds for parallel first-order stochastic methods in the smooth nonconvex regime, thereby achieving optimal time complexity under data heterogeneity and without restrictive similarity assumptions.
    Our analysis further establishes that \algn remains optimal under arbitrary and even time-varying worker computation speeds, closing a fundamental gap in the theory of asynchronous optimization.
\end{abstract}
\section{Introduction}
Modern machine learning increasingly depends on large-scale distributed training across clusters with hundreds or even thousands of GPUs \citep{shoeybi2019megatron, GPT3, narayanan2021efficient}.
However, classical synchronous training methods struggle to scale in these settings, as device failures, network instabilities, and synchronization overheads introduce significant inefficiencies \citep{chen2016revisiting, grattafiori2024llama}.
These issues become even more pronounced in environments with heterogeneous computational power, such as Federated Learning (FL), where devices range from high-end datacenter GPUs to resource-constrained edge hardware \citep{konevcny2016federated, mcmahan2016federated, li2020federated, kairouz2021advances}.
Because synchronous methods are bottlenecked by the slowest participants, faster devices remain idle, leading to severe underutilization of computational resources when stragglers---nodes slowed down by computation or communication---lag significantly behind.

One way to reduce synchronization bottlenecks is to equip data centers with homogeneous GPUs.  
However, this approach is prohibitively expensive and difficult to scale: upgrading to faster GPUs would require replacing all devices simultaneously, since heterogeneous hardware cannot be combined efficiently.  
Even then, homogeneity does not eliminate synchronization issues, as hardware failures and device dropouts during training still cause stragglers and idle time.  
Moreover, this solution applies only to controlled datacenter environments and is infeasible in FL, where edge devices are outside the server's control.

A more promising approach is to shift from hardware solutions to algorithmic ones by adopting asynchronous optimization methods.
These methods remove the need for synchronization, allowing fast workers to contribute updates without waiting for slower ones \citep{tsitsiklis1986distributed, recht2011hogwild, agarwal2011distributed, dean2012large, li2014communication}.
Despite their appeal, asynchronous methods are more difficult to analyze.
In particular, a meaningful analysis would require studying \textit{time to convergence}, rather than iteration complexity only.
While iteration complexity is the traditional metric in optimization, it does not necessarily reflect real-world training speed in parallel settings: a method that performs more iterations may finish faster in wall-clock time if those iterations can be computed without waiting for slow workers.
This distinction raises a fundamental question: \textit{among all parallel methods, which ones are provably fastest in theory?}
To make this question precise, we restrict our attention to smooth nonconvex problems and to stochastic first-order methods, encompassing algorithms with or without synchronization.
This will be the only setting considered in this paper.

\begin{table*}
    \caption{
        Comparison of time complexities for parallel first-order methods under the fixed computation time model, where each worker $i$ takes a fixed time $\tau_i$ to compute a stochastic gradient, with the times ordered so that $\tau_n$ is the largest \eqref{eq:fixed_time}.
        We denote by $\tauavg \eqdef \tfrac{1}{n}\sum_{i=1}^n \tau_i$ the average computation time across all workers.
        The table shows how the time complexity of different algorithms depends on key problem parameters: the initial function suboptimality $\Delta \eqdef f(x^0) - f^*$ (\Cref{ass:lower_bound}), the target stationarity $\varepsilon$, the variance bound of the stochastic gradients $\sigma^2$ (\Cref{ass:stochastic_variance_bounded}), and smoothness constants.
        Specifically, $L_f$ is the smoothness constant of $f$ (\Cref{eq:lipschitz_gradient}); $L_{\max} \eqdef \max_{i\in[n]} L_{f_i}$ with $L_{f_i}$ the smoothness constant of $f_i$; and $L$ is a constant associated with our new smoothness-type assumption (\Cref{ass:lipschitz_constant}).
        They satisfy $L_f \le L \le L_{\max}$ (\Cref{lemma:smoothness_relation}).
        All stated time complexities hide universal constant factors.\\
        Each column indicates whether a method satisfies the following desirable properties:
        \textbf{Optimal:} achieves the theoretical lower bound derived by \citet{tyurin2024optimal} for parallel first-order stochastic methods in heterogeneous data setting.
        \textbf{No sync.:} does not require synchronization and is therefore \textit{asynchronous}.
        \textbf{No idle workers:} all workers remain busy without waiting, so computational resources are fully utilized.
        \textbf{No discarded work:} no computation is wasted, and no worker is stopped mid-computation.
        Our new method, \algn, is the first asynchronous method to achieve optimal time complexity, while also ensuring full resource utilization (no idle workers) and no discarded computations / work.
    }
    \label{table:fixedtime}
    \centering
    \begin{threeparttable}
      \resizebox{\textwidth}{!}{
        \begin{tabular}{cccccc}
            \toprule
                \bf Method 
                & \bf Time Complexity 
                & \bf Optimal 
                & \bf No sync.
                & \bf No idle workers
                & \bf No discarded work\\
            \midrule
                \makecell{\naiveminibatch \\ (\Cref{sec:minibatch})}
                & $\frac{L_f\Delta}{\varepsilon}\(\tau_n + \tau_n \frac{\sigma^2}{n\varepsilon}\)$
                & \crossmarkred
                & \crossmarkred
                & \crossmarkred
                & \checkmarkgreen \\
            \midrule
                \makecell{\iasgd \\ \citep{wang2025incremental} \\ (\Cref{sec:iasgd_time_complexity})}
                & $\frac{L_{\max}\Delta}{\varepsilon}\(\tau_n + \tau_n \frac{\sigma^2}{n\varepsilon}\)$ \textsuperscript{{\color{mydarkblue} (\dag)}}
                & \crossmarkred
                & \checkmarkgreen
                & \checkmarkgreen
                & \checkmarkgreen \\
            \midrule
                \makecell{\malenia \\ \citep{tyurin2024optimal}}
                & $\frac{L_f\Delta}{\varepsilon}\(\tau_n + \tauavg \frac{\sigma^2}{n\varepsilon}\)$
                & \checkmarkgreen
                & \crossmarkred
                & \checkmarkgreen
                & \crossmarkred \\
            \midrule
                \makecell{\algn \textbf{(new)} \\ (\Cref{algo:Ringleader}; \Cref{thm:time_complexity})}
                & $\frac{L\Delta}{\varepsilon}\(\tau_n + \tauavg \frac{\sigma^2}{n\varepsilon}\)$
                & \checkmarkgreen \textsuperscript{{\color{mydarkblue} (\ddag)}}
                & \checkmarkgreen
                & \checkmarkgreen
                & \checkmarkgreen\\
                \midrule
                \midrule
                    \makecell{Lower Bound \\ \citep{tyurin2024optimal}}
                    & $\frac{L_f\Delta}{\varepsilon}\(\tau_n + \tauavg \frac{\sigma^2}{n\varepsilon}\)$
                    & \textbf{---}
                    & \textbf{---}
                    & \textbf{---}
                    & \textbf{---}\\
                \bottomrule
            \end{tabular}%
        }

    \resizebox{\textwidth}{!}{
        \begin{minipage}{\textwidth}
            \begin{tablenotes}[para,flushleft]
                \footnotesize
                \item[{\color{mydarkblue} (\dag)}]
                    The analysis presented by \citet{wang2025incremental} is carried out under the assumption that each $f_i$ is smooth with the same smoothness constant, which is equivalent to requiring that each $f_i$ is $L_{f_i}$--smooth and then using the upper bound $L_{\max}$ for all $L_{f_i}$.
                    However, this strict assumption is not necessary: they could instead use our relaxed smoothness-type assumption (\Cref{ass:lipschitz_constant}), in which case the constant improves to $L$, and their analysis remains unchanged.
                \\
                \item[{\color{mydarkblue} (\ddag)}]
                    The time complexities of \algn and \malenia differ in the smoothness constant only.
                    Since \malenia is optimal, \algn is also optimal whenever $L$ exceeds $L_f$ by at most a universal constant factor, that is, $L = \cO(L_f)$.
            \end{tablenotes}
        \end{minipage}
    }
    \end{threeparttable}
\end{table*} 

Recently, \citet{tyurin2024optimal} studied this very regime, where they derived lower bounds.
They then proposed two algorithms: \rennala, designed for the \textit{homogeneous data setting}, where all workers draw samples from the same distribution, and \malenia, for the \textit{heterogeneous data setting}, where data distributions differ across workers.
They showed that both methods are optimal---achieving the lower bounds---and, perhaps surprisingly, both are synchronous (they periodically synchronize the workers).
The key idea in both is to fully utilize the available computational resources by keeping workers continuously busy: each worker computes independently, and synchronization occurs only after a sufficient number of gradient computations have been accumulated.

At first, the result of \citet{tyurin2024optimal} suggested a rather pessimistic outlook for asynchronous methods: despite their practical appeal, they showed that existing asynchronous methods are not optimal and that the method achieving the lower bound is \textit{synchronous}.
This created the view that optimality is inherently tied to synchronization.
However, this view was overturned by \citet{maranjyan2025ringmaster}, who, in the \textit{homogeneous data setting}, introduced \ringmaster---the first asynchronous \sgd method to achieve the same optimal time complexity as the synchronous \rennala.
Although both methods share the same theoretical guarantees, \ringmaster can be faster than \rennala in practice, since it avoids synchronization and benefits from more frequent updates.

Nevertheless, the work of \citet{maranjyan2025ringmaster} established optimality in the \textit{homogeneous data setting} only.
The question of whether some variant of a parallel method that does not rely on synchronization (i.e., is asynchronous) can also be optimal in the more general \textit{heterogeneous data setting} remained open.
In this work, we close this gap and answer the question affirmatively.

The heterogeneous data setting is both important and practically relevant.  
In FL, for instance, such heterogeneity arises naturally as participants hold distinct datasets \citep{zhao2018federated, li2020federated, tan2022towards}.
Yet this setting is significantly more challenging than the homogeneous one.
The standard philosophy of asynchronous \sgd---updating the model after every gradient computation---can be harmful here: fast workers contribute updates more frequently, causing the optimization process to become biased toward their local data.
To mitigate this, most existing asynchronous methods address this issue by assuming similarity across client distributions \citep{mishchenko2022asynchronous, koloskova2022sharper, nguyen2022federated, islamov2024asgrad}.
While this assumption simplifies the analysis, it is often unrealistic in practice, where clients may involve fundamentally different populations (e.g., hospitals with distinct demographics, mobile users in different countries, or financial institutions under varied regulations).

Recent work by \citet{wang2025incremental} took an important step toward removing these restrictive assumptions by proposing Incremental Aggregated Asynchronous \algname{SGD} (\iasgd), a method that provably converges without similarity assumptions.
However, their method achieves the same time complexity as standard \naiveminibatch (see the first two rows of \Cref{table:fixedtime})---the simplest synchronous \sgd baseline, which waits to collect one gradient from each worker before every update---thus failing to provide the computational advantages that motivate asynchronous approaches in the first place.

To the best of our knowledge, the only method proven to be optimal in the \textit{heterogeneous data setting} is the synchronous algorithm \malenia of \citet{tyurin2024optimal}, which, notably, does not rely on similarity assumptions.
However, synchronization is a major bottleneck in practice: although synchronous and asynchronous methods can share the same theoretical complexity, asynchronous methods are often faster in practice because they avoid costly synchronization and benefit from more frequent updates, as demonstrated in the homogeneous case by \citet{maranjyan2025ringmaster}.

This raises a fundamental question: \textit{Is it possible to design an asynchronous method that requires no similarity assumptions while still achieving optimal time complexity?}
In this paper, we answer this question affirmatively by introducing \algn, the first asynchronous \sgd method that achieves optimal time complexity\footnote{
    Throughout the paper, we refer to our method as optimal.
    Formally, this holds whenever the constant $L$---associated with our new smoothness-type assumption (\Cref{ass:lipschitz_constant})---is at most a constant factor larger than the smoothness constant $L_f$ used in the derived lower bounds \citep{tyurin2024optimal}.
    See \Cref{table:fixedtime} for details.
    }
in the \textit{heterogeneous data setting}.
Importantly, \algn attains this without relying on restrictive similarity assumptions.
\subsection{Contributions}
Our main contributions are the following:
\begin{itemize}
    \item \textbf{Optimal asynchronous SGD under data heterogeneity.}
    We prove that \algn (\Cref{algo:Ringleader}) is, to the best of our knowledge, the first asynchronous method in the heterogeneous setting under the fixed computation model \eqref{eq:fixed_time} matching the lower bounds for parallel methods of \citet{tyurin2024optimal}, when the smoothness-type constant~$L$ in \Cref{ass:lipschitz_constant} is within a constant factor of the smoothness~$L_f$ used to obtain the lower bounds (\Cref{table:fixedtime}).
    Importantly, \algn attains this without any similarity assumptions across clients' data.
    \item \textbf{Additional useful properties.}  
    Beyond achieving optimal time complexity, our method \algn satisfies two additional desirable properties: (i) all workers remain continuously active (\textit{no idle workers}), and (ii) every computed gradient is incorporated into the update (\textit{no discarded work}).
    These properties are crucial in practice, as they ensure maximum resource utilization: all workers contribute at all times, and their computations are never wasted.
    \Cref{table:fixedtime} compares \algn against benchmark algorithms with respect to these properties.
    \item \textbf{Parameter-free design.}
    In contrast to the optimal synchronous method \malenia \citep{tyurin2024optimal}, which requires prior knowledge of the gradient variance bound and target accuracy, our method operates in the fixed computation time model without such parameters (except for the stepsize, needed only to match the optimal theoretical rate).
    This makes it far more practical for real-world deployments, where these quantities are typically unknown or difficult to estimate.
    The same parameter-free improvement can also be extended to \malenia, as we discuss in \Cref{sec:malenia_param_free}.
    \item \textbf{Universal computation model.}
    In \Cref{sec:arbitrary_time}, we extend our analysis beyond the fixed computation time model to the general setting of arbitrarily varying computation times, accommodating virtually any computational behavior, including stochastic or adversarial patterns, while retaining optimality time complexity.
    \item \textbf{Empirical validation.}
    In \Cref{sec:experiments}, we evaluate \algn against benchmark methods on illustrative toy problems.
    The results validate our theoretical findings and demonstrate clear practical advantages over the baselines.
\end{itemize}
\section{Problem Setup}
We consider a distributed learning setting with $n$ workers, where each worker~$i$ possesses its own local data distribution~$\cD_i$.
Our goal is to solve the following distributed optimization problem:
\begin{equation}\label{eq:problem}
    \minimize\limits_{x \in \R^d}
    \left\{ f(x) \eqdef \frac{1}{n}\sum_{i=1}^{n}f_i(x) \right\},
    \quad\text{where}\quad
    f_i(x) \eqdef \ExpSub{\xi_i \sim \cD_i}{f_i(x; \xi_i)}.
\end{equation}
Here $f_i \colon \R^d \to \R$ denotes the local objective of worker $i$, defined as the expectation of the sample loss $f_i(x;\xi_i)$ over data points $\xi_i$ drawn from its local distribution $\cD_i$.
\subsection{Worker Heterogeneity Model}\label{sec:compute_times}
We first focus on the case where workers have constant computation speeds, as this setting is more intuitive and serves as a foundational model for understanding the dynamics of asynchronous distributed optimization.
The extension to arbitrary computation times is presented in \Cref{sec:arbitrary_time}.

Following the \textit{fixed computation model} \citep{mishchenko2022asynchronous}, we formalize:
\begin{equation}
    \label{eq:fixed_time}
    \begin{aligned}
        &\text{Each worker } i \text{ requires } \tau_i \text{ seconds\footnotemark\ to compute one stochastic gradient} \ \nabla f_i(x,\xi_i)~.\\
        &\text{Without loss of generality, assume } 0 < \tau_1 \leq \tau_2 \leq \cdots \leq \tau_n~.
    \end{aligned}
\end{equation}
\footnotetext{
    One could alternatively assume that each worker requires \textit{at most} $\tau_i$ seconds. 
    Under this formulation, all of our upper bounds would still hold; however, the lower bound would no longer be valid. 
    For this reason, we adopt the assumption that each worker requires exactly $\tau_i$ seconds.}
We assume communication is infinitely fast (taking $0$ seconds), both from workers to the server and from the server to workers\footnote{
    Alternatively, one could define $\tau_i$ as the time required for a worker to both compute a gradient and communicate it to the server, while keeping server-to-worker communication infinitely fast.
    Our upper bounds would still hold under this formulation, but the lower bounds would no longer apply, so we use the simpler model.}.
This is a modeling choice---arguably the simplest one---and has been the standard assumption in prior work \citep{mishchenko2022asynchronous, koloskova2022sharper, tyurin2024optimal, maranjyan2025ringmaster}, even if not always stated explicitly.
We further discuss the motivation and limitations of this abstraction in \Cref{sec:communication_limitations}.
A related study by \citet{tyurin2024shadowheart} considers the case where communication is non-negligible and proposes techniques to address it.

Finally, we denote by $\tauavg \eqdef \tfrac{1}{n}\sum_{i=1}^n \tau_i$ the average computation time across all workers.
\subsection{Notations}
We denote the standard inner product in $\R^d$ by 
$$
    \inp{x}{y} = \sum_{i=1}^{d} x_i y_i~,
$$ 
and the corresponding Euclidean norm by $\norm{x} \eqdef \sqrt{\inp{x}{x}}$.
We use $[n] \eqdef \{1,2,\dots,n\}$ to denote the index set, and $\Exp{\cdot}$ for mathematical expectation.
For functions $\phi, \psi : \mathcal{Z} \to \R$, we write $\phi = \cO(\psi)$ if there exists a constant $C > 0$ such that $\phi(z) \leq C \psi(z)$ for all $z \in \mathcal{Z}$.
\subsection{Assumptions}
\label{sec:assumption}
We consider the following standard assumptions for the nonconvex setting.
\begin{boxedassumption}
    \label{ass:stochastic_variance_bounded}
    For each $i \in [n]$ and every $\xi$, the function $f_i(x;\xi)$ is differentiable with respect to its first argument $x$.
    Moreover, the stochastic gradients are unbiased and have bounded variance $\sigma^2 \geq 0$, that is,
    \begin{gather*}
         \ExpSub{\xi_i \sim \cD_i}{\nabla f_i(x;\xi_i)} = \nabla f_i(x), 
         \quad \forall x \in \R^d, \;\; \forall i \in [n]~,\\
         \ExpSub{\xi_i \sim \cD_i}{\sqnorm{\nabla f_i(x;\xi_i) - \nabla f_i(x)}} \leq \sigma^2,\quad \forall x \in \R^d, \;\; \forall i \in [n]~.
    \end{gather*}
\end{boxedassumption}
\begin{boxedassumption}\label{ass:lipschitz_constant}
    Each function $f_i$ is differentiable.
    There exists a constant $L>0$ such that, for all $x \in \R^d$ and $y_1,\dots,y_n \in \R^d$, 
    $$
        \sqnorm{\nabla f(x) - \frac{1}{n} \sum_{i=1}^n \nabla f_i(y_i)}
        \le \frac{L^2 }{n} \sum_{i=1}^n \sqnorm{x - y_i}.
    $$
\end{boxedassumption}
Recall the standard definition of smoothness
\begin{boxeddefinition}[Smoothness]\label{eq:lipschitz_gradient}
    A differentiable function $\phi \colon \R^d \to \R$ is called $L_\phi$--smooth if
    \begin{equation*}
        \norm{\nabla \phi(x) - \nabla \phi(y)} \le L_\phi \norm{x - y}, \quad \forall x,y \in   \R^d~.
    \end{equation*}
    By convention, $L_\phi$ denotes the smallest such constant.
\end{boxeddefinition}
Note that \Cref{ass:lipschitz_constant} is stronger than requiring $f$ itself to be $L_f$--smooth, yet weaker than all $f_i$ being $L_{f_i}$--smooth.
The following lemma establishes the relation among the constants $L_f$, $L$, and $L_{f_i}$.
\begin{boxedlemma}[Smoothness Bounds; Proof in \Cref{proof:smoothness}]\label{lemma:smoothness_relation}
    Suppose \Cref{ass:lipschitz_constant} holds with constant $L>0$.
    Then $f$ is $L_f$--smooth with $L_f \le L$.
    Moreover, if each $f_i$ is $L_{f_i}$--smooth, then \Cref{ass:lipschitz_constant} is satisfied, and we have
    $$
        L_f \;\le\; L \;\le\; \sqrt{\frac{1}{n}\sum_{i=1}^n L_{f_i}^2} \;\le\; \max_{i \in [n]} L_{f_i} =: L_{\max}~.
    $$
    Finally, if all $f_i$ are identical, i.e., $f_i = f$ for all $i \in [n]$, then $L = L_f$.
\end{boxedlemma}
To the best of our knowledge, prior work on asynchronous \sgd under data heterogeneity has always assumed smoothness of each $f_i$, including works by \citet{koloskova2022sharper,nguyen2022federated,wang2025incremental}.
\begin{boxedassumption}\label{ass:lower_bound}
    There exists $f^* > -\infty$ such that $f(x) \geq f^*$ for all $x \in \R^d$.
    We define $\Delta \eqdef f(x^0) - f^*,$ where $x^0$ is the starting point of the optimization methods.
\end{boxedassumption}
Under these assumptions; our objective is to find an $\varepsilon$--stationary point---a (possibly random) vector $x$ satisfying $\E{\sqnorm{\nabla f(x)}} \leq \varepsilon$.
\section{Background and Motivation}
In this section, we review relevant existing methods in distributed optimization and discuss their limitations to motivate our algorithmic design.
We begin with \naiveminibatch as the canonical synchronous baseline, then consider \malenia \citep{tyurin2024optimal}, the first synchronous \sgd method to attain optimal time complexity.
We then turn to the challenges of asynchronous approaches, focusing on the recent \iasgd \citep{wang2025incremental} and its limitations.
Finally, we outline how these limitations can be addressed and introduce the core idea behind our algorithm.
\subsection{\naiveminibatchtitle}
\label{sec:minibatch}
\naiveminibatch provides the most straightforward approach to solving problem \eqref{eq:problem} in a distributed setting.
\paragraph{Algorithm Description.}
At each iteration $k$, the algorithm performs the following update:
$$
    x^{k+1} = x^k - \gamma \frac{1}{n} \sum_{i=1}^n \nabla f_i \(x^k; \xi_i^k\).
$$
The algorithm operates synchronously: at each iteration, the server waits to receive one stochastic gradient from each of the $n$ workers, all of which are evaluated at the current iterate $x^k$.
Once all gradients are collected, the server constructs an unbiased minibatch estimator of the full gradient by averaging these stochastic gradients and performs a standard gradient descent step.
\paragraph{Limitations.}
This synchronous approach has a significant computational bottleneck.
Each iteration requires waiting for the slowest worker to complete its gradient computation, resulting in the iteration time $\tau_n \eqdef \max_{i \in [n]} \tau_i$ \eqref{eq:fixed_time}.
Consequently, faster workers remain idle while waiting for stragglers, which leads to inefficient utilization of available computational resources.
\paragraph{Theoretical Performance.} 
From a convergence perspective, \naiveminibatch achieves the iteration complexity
$
    \cO\( \nicefrac{L_f\Delta}{\varepsilon}\(1 + \nicefrac{\sigma^2}{n\varepsilon}\)\)
$
to reach an $\varepsilon$--stationary point \citep{cotter2011better,goyal2017accurate,gower2019sgd}.
The corresponding time complexity becomes
$$
    \cO\(\frac{\tau_n L_f\Delta}{\varepsilon}\(1 + \frac{\sigma^2}{n\varepsilon}\)\).
$$
This motivates the development of methods that can better utilize fast workers without waiting for stragglers.
\subsection{\maleniatitle}
\malenia \citep{tyurin2024optimal} addresses the straggler problem of \naiveminibatch by ensuring continuous worker utilization, rather than waiting for the slowest worker in each iteration.
However, it is fundamentally a \minibatch algorithm: in every iteration, it constructs an unbiased gradient estimator from multiple gradients and then performs a synchronous update.
The key distinction lies in how the batch is collected---\malenia gathers gradients asynchronously, allowing potentially multiple contributions from the same worker within a single iteration, unlike \naiveminibatch.
\paragraph{Algorithm Description.}
At each iteration $k$, all workers continuously compute stochastic gradients at the current point $x^k$.
The server accumulates these gradients until the stopping condition
\begin{equation}\label{eq:malenia_condition}
    \(\frac{1}{n} \sum_{i=1}^n \frac{1}{b_i^k} \)^{-1} \ge \max\left\{1,\frac{\sigma^2}{n\varepsilon}\right\}
\end{equation}
is satisfied, where $b_i^k$ denotes the number of gradients received from worker $i$.
Once this condition is met, the server performs the update
$$
    x^{k+1}
    = x^k - \gamma \bar g^k
    \eqdef x^k - \gamma \frac{1}{n} \sum_{i=1}^n \bar g_i^k~,
$$
where $\bar g_i^k$ is the average of the stochastic gradients received from worker $i$
$$
    \bar g_i^k = \frac{1}{b_i^k} \sum_{j=1}^{b_i^k} \nabla f_{i}\(x^{k}; \xi_{i}^{k,j}\).
$$
Because stochastic gradients are collected asynchronously, but the update is performed only after all required gradients are received, the algorithm can be regarded as semi-asynchronous---asynchronous in gradient collection but synchronous in parameter updates.
\paragraph{Stopping Condition Rationale.}
The left-hand side of condition \eqref{eq:malenia_condition} appears in the variance of the gradient estimator $\bar g^k$.
Ensuring that this quantity is sufficiently large allows the algorithm to control the variance at each step.
Moreover, condition \eqref{eq:malenia_condition} guarantees that every worker computes at least one gradient, which in turn yields a smaller variance than that of the estimator in \naiveminibatch.
\paragraph{Theoretical Performance.}
This asynchronous gradient collection strategy achieves the optimal time complexity
$$
    \cO\(\frac{L_f\Delta}{\varepsilon}\(\tau_n + \tauavg \frac{\sigma^2}{n\varepsilon}\)\),
$$
as shown by \citet{tyurin2024optimal}.
The key improvement over \naiveminibatch is that the variance term $\sigma^2$ is now multiplied by $\tau_{\text{avg}}$ instead of $\tau_n$.
Since $\tau_{\text{avg}} \ll \tau_n$ in computationally heterogeneous environments, \malenia can potentially provide substantial speedup in highly heterogeneous regimes.

The algorithm's benefits are most pronounced in the high-noise settings (where $\nicefrac{\sigma^2}{n \varepsilon}$ is large).
In low-noise scenarios or when $\sigma = 0$, \malenia offers no advantage over \naiveminibatch, since collecting multiple gradients per worker provides no benefit in terms of variance reduction.
\paragraph{Limitations.}
The main limitation of \malenia is its synchronous nature.
After collecting the necessary gradients, the server must broadcast the updated model to all workers simultaneously.
This all-at-once synchronization creates substantial communication overhead, which can become a major scalability bottleneck in large-scale or bandwidth-limited environments.

Moreover, the synchronization process forces the server to discard ongoing gradient computations from workers who are actively computing when the broadcast occurs.
Since workers operate continuously during the gradient collection phase, they must abandon their current computations upon receiving the new model, wasting valuable computational resources that could otherwise contribute to convergence.

Additionally, \malenia requires prior knowledge of the noise level $\sigma$ and the target accuracy $\varepsilon$ to evaluate the stopping condition \eqref{eq:malenia_condition}.
This dependence on problem-specific parameters, which are often unknown or difficult to estimate in practice, significantly limits its practical applicability.

These synchronization bottlenecks motivate the development of asynchronous optimal methods.
Beyond avoiding coordination overhead, asynchronous approaches enable immediate model updates upon gradient arrival, which can accelerate convergence through more frequent parameter updates.
This immediate processing is particularly beneficial for sparse models where gradients typically affect disjoint parameter subsets, allowing parallel updates without interference \citep{recht2011hogwild}.
\subsection{Toward Asynchronous Methods}
A naive approach to making optimization asynchronous would be to update the model immediately upon receiving any gradient, using
$$
    x^{k+1} = x^k - \gamma \nabla f_{i_k}\(x^{k-\delta^k};\xi_{i_k}^{k-\delta^k}\),
$$
where $i_k$ denotes the worker that sent the gradient at iteration $k$, and $\delta^k \ge 0$ is its delay, i.e., the number of server updates that occurred while the gradient was being computed.
Delays arise naturally in asynchronous execution: fast workers return gradients quickly and proceed with updated models, while slower workers compute on stale iterates; by the time they return their gradients, several server updates may already have occurred.
Consequently, $\delta^k$ is only determined when the server receives a gradient: at that point, it knows both the model iterate used by the worker and the current server iterate, and $\delta^k$ is simply the difference between the two.

\paragraph{Limitations.}
This approach suffers from a fundamental bias problem when workers have heterogeneous data distributions.
Faster workers send updates more frequently, causing their local objectives to dominate the optimization and pull the model toward their own minimizers.
Slower workers, when their updates finally arrive, push the model back toward different solutions.
As a result, the iterates may oscillate or stagnate, failing to converge to a stationary point.

This bias also makes theoretical analysis difficult.
Classical \algname{SGD}-style proofs rely on one-step progress toward minimizing the global function, but here each update direction reflects a different local objective.
Without additional data similarity assumptions \citep{mishchenko2022asynchronous,koloskova2022sharper,nguyen2022federated,islamov2024asgrad}, it becomes impossible to extend the analysis to the global function---yet such assumptions are rarely realistic when data can be arbitrarily heterogeneous across machines or organizations.

The root cause is that each update uses a gradient from one worker only.
A better strategy is to incorporate information from all workers, even if some gradients are computed at stale iterates.
This idea motivates methods such as Incremental Aggregated Gradient (\algname{IAG}) \citep{blatt2007convergent,gurbuzbalaban2017convergence,vanli2018global} and Stochastic Averaged Gradient (\algname{SAG}) \citep{roux2012stochastic,schmidt2017minimizing}, which maintain and aggregate gradients from all workers.






\subsection{\iasgdtitle}
As discussed above, the key insight is to construct a gradient estimator using information from all workers at each model update.
\iasgd \citep{wang2025incremental} achieves this by maintaining a gradient table on the server, similar to \algname{SAG} or \algname{IAG}, but with asynchronous table updates.
\paragraph{Algorithm Overview.} 
The server maintains a gradient table $\{g_i\}_{i=1}^n$ that stores the most recent gradient received from each of the $n$ workers.
The table is initialized by collecting one gradient from each worker at the starting point $x^0$.
The server then computes the first model update 
$$
    x^1 = x^0 - \gamma \bar g \eqdef x^0 - \gamma \frac{1}{n} \sum_{i=1}^n g_i
$$
and broadcasts $x^1$ to all workers.
Subsequently, the workers compute stochastic gradients in parallel, and their corresponding table entries are updated asynchronously as soon as the computations finish.

At each iteration $k$, the server performs the following steps:
\begin{enumerate}
    \item Receive gradient 
        $\nabla f_{i_k}\(x^{k-\delta_{i_k}^k};\xi_{i_k}^{k-\delta_{i_k}^k}\)$ 
        from worker $i_k$
    \item Update the gradient table entry: 
        $g_{i_k} = \nabla f_{i_k}\(x^{k-\delta_{i_k}^k};\xi_{i_k}^{k-\delta_{i_k}^k}\)$
    \item Perform model update: 
        $x^{k+1} = x^k - \gamma \bar g = x^k - \gamma \frac{1}{n} \sum_{i=1}^n g_i$
    \item Send the updated iterate $x^{k+1}$ to worker $i_k$ for its next computation
\end{enumerate}
The gradient estimator $\bar g$ combines the most recent gradient from each worker, ensuring that every update reflects information from the entire set of workers despite asynchronous execution.
In this way, the method retains the statistical benefits of using all workers' data while allowing them to operate independently, thereby avoiding the synchronization bottlenecks that limit scalability.

Note that, due to asynchrony, the stochastic gradients stored in the table generally correspond to different iterates of the model.
We therefore record each worker $i$'s delay $\delta_i^k$ to track the iterate at which its gradient was computed.

\paragraph{Theoretical Performance.}
The iteration complexity of this algorithm was established by \citet{wang2025incremental}, and by a straightforward conversion to the fixed computation time model (see \Cref{sec:iasgd_time_complexity}), we obtain the corresponding time complexity
$$
    \cO\(\frac{\tau_nL_{\max}\Delta}{\varepsilon}\(1 + \frac{\sigma^2}{n\varepsilon}\)\),
$$
which matches the complexity of \naiveminibatch.
This indicates that asynchronous execution alone does not provide computational advantages over synchronous approaches.
Thus, a fundamental challenge lies in how gradient delay affects convergence, which we address next.
\subsection{Motivation}
The primary limitation of asynchronous algorithms stems from gradient delay, which can significantly degrade convergence performance.
Large delays can cause the optimization steps to follow suboptimal trajectories, which disrupts convergence.

This delay problem occurs even in homogeneous settings where all functions $f_i$ are equal ($f_i \equiv f$ for all $i\in [n]$).
The state-of-the-art solution for this case, \algname{Ringmaster ASGD} \citep{maranjyan2025ringmaster}, achieves optimal time complexity by explicitly controlling delay to prevent it from becoming large.
\algname{Ringmaster ASGD} accomplishes this by discarding gradients that arrive with large delays.

Unfortunately, this gradient-discarding strategy fails in the heterogeneous setting of \iasgd.
The fundamental issue is that slow workers inevitably experience large delays due to their computational constraints.
If we ignore or discard their delayed gradients, the corresponding table entries remain outdated and may never be updated again if subsequent gradients also arrive late and are discarded.
This creates a persistent information bottleneck that degrades the quality of the gradient estimator and harms convergence.

This issue suggests we should prevent such situations from occurring by controlling the number of updates performed using fast workers.
The simplest approach would be to ignore some updates from fast workers, but this contradicts the core principle of asynchronous methods whereby all workers compute gradients continuously.

Instead, our approach \textit{buffers} the gradients produced by fast workers rather than applying them immediately, similar to the strategy in \malenia.
By buffering gradients and performing a model update only once a sufficient number have been collected, we control the delays induced by slow workers while keeping all workers continuously active.
This buffering mechanism provides an additional advantage: when multiple gradients computed at the same iterate are aggregated, they yield lower-variance estimates, thereby improving convergence.
\section{\red Ringleader ASGD}
We now formally introduce our algorithm, \algn.

\algn builds upon three key insights from existing methods.
First, inspired by \iasgd, we maintain a gradient table\footnote{
    It is not strictly necessary to maintain a table on the server.
    An alternative, as in \iasgd \citep{wang2025incremental}, is to store one vector on each worker along with an additional vector on the server.
    However, this can be problematic when workers have limited memory.
    For clarity and simplicity, we adopt the server-side table formulation in our description.}
to ensure that information from all workers is incorporated into every update, which eliminates the need for data similarity assumptions between workers.
Second, following \ringmaster, we recognize that controlling gradient delay is essential for efficient asynchronous optimization.
Third, drawing from \malenia, we use buffering of stochastic gradients---rather than discarding delayed ones---to control delays while preserving valuable computations, enabling continuous utilization of all resources.

An important property of the algorithm is that all workers remain continuously active, computing stochastic gradients.
As soon as a worker finishes computing a gradient, it immediately sends it to the server.
Recall that we assumed communication is instantaneous, i.e., takes zero time (\Cref{sec:compute_times}).
When the server receives a gradient, it either buffers it for later use or applies it immediately to perform a model update.
If the gradient is buffered, no further action is taken and the worker simply continues computing and sending new gradients.
If the server decides to perform an update, it updates the model and sends the updated model back to the worker that provided the gradient.
This server-to-worker communication is also assumed instantaneous, after which the worker resumes computing gradients at the new model point, ensuring that workers are never idle.

The algorithm proceeds in rounds.
In each round, exactly $n$ model updates are performed---one for each worker.
Specifically, when a worker sends a stochastic gradient, the server may apply an update and return the updated model to that worker, but it ensures that each worker receives an updated model at most once per round.
Repeating this procedure $n$ times ensures that every worker obtains exactly one fresh model per round, which in turn keeps delays bounded.
To avoid discarding the computations of fast workers, the server buffers their gradients and applies them only at the appropriate moment, thereby guaranteeing that each round consists of exactly $n$ updates.

Each round consists of two phases:
\begin{itemize}
    \item \textbf{Phase~1:} Buffer stochastic gradients in a table until at least one gradient from each worker is available.
    \item \textbf{Phase~2:} Perform exactly $n$ updates (one for each worker), then discard the old stochastic gradients from the table and return to Phase~1 to start collecting fresh ones.
\end{itemize}
The complete algorithm is formalized in \Cref{algo:Ringleader}.

\begin{figure}[t]
\centering
\begin{minipage}{\linewidth}
\begin{algorithm}[H]\caption{\algn (server algorithm)}
    \label{algo:Ringleader}
    \begin{algorithmic}[1]
        \STATE \textbf{Input:} Stepsize $\gamma > 0$, initial point $x^0 \in \R^d$
        \STATE \textbf{Initialization:} Broadcast $x^0$ to all workers, which then start running \Cref{algo:Ringleader_worker} in parallel
        \STATE Set $k = 0$, $S = \emptyset$; initialize $G_i = 0$, $b_i = 0$ for all $i \in [n]$
        \WHILE{True}
        \STATE {\orange \textbf{--- Phase~1: await stochastic gradients from all workers ---}}
        \WHILE{$S \ne [n]$}
            \STATE Receive stochastic gradient $g_j^k$ (computed at $x^{k-\delta_j^k}$) from some worker $j \in [n]$
                \label{line:receive_ph1}
            \STATE $G_j = G_j + g_j^k$;\;  $b_j = b_j + 1$;\; $S = S \cup \{j\}$
                \label{line:update_ph1}
        \ENDWHILE
        \STATE {\orange \textbf{--- Phase~2: perform exactly one update for every worker ---}}
        \STATE $x^{k+1} = x^k - \gamma \frac{1}{n}\sum_{i=1}^n \nicefrac{G_i}{b_i}$
            \hfill $\diamond$ Update using averaged gradients from all workers
            \label{line:step_last_worker}
        \STATE Broadcast $x^{k+1}$ to worker $j$
            \hfill $\diamond$ Last worker to complete Phase~1
            \label{line:broadcast_last_worker}
        \STATE $k=k+1$;\; $S = S \setminus \{j\}$
            \label{line:counter_last_worker}
        \STATE $g_i^+ = 0$, $b_i^+ = 0$ for all $i \in [n]$;\; $S^+ = \emptyset$
            \hfill $\diamond$ Initialize temporary buffer for the next round
            \label{line:buffer_initialization}
        \WHILE{$S \ne \emptyset$}
            \STATE Receive stochastic gradient $g_j^k$ (computed at $x^{k-\delta_j^k}$) from some worker $j \in [n]$
            \IF{$j \in S$}
                \STATE $G_j = G_j + g_j^k$;\; $b_j = b_j + 1$
                    \label{line:update_table_ph2}
                \STATE $x^{k+1} = x^k - \gamma \frac{1}{n}\sum_{i=1}^n \nicefrac{G_i}{b_i}$
                    \label{line:step_ph2}
                \STATE Broadcast $x^{k+1}$ to worker $j$
                    \label{line:broadcast_ph2}
                \STATE $k=k+1$;\; $S = S \setminus \{j\}$
                    \label{line:counter_ph2}
            \ELSE
                \STATE $G_j^+ = G_j^+ + g_j^k$;\; $b_j^+ = b_j^+ + 1$;\; $S^+ = S^+ \cup \{j\}$ 
                    \hfill $\diamond$ Buffer for next round
                    \label{line:buffer}
            \ENDIF
        \ENDWHILE
        \STATE $G_i = G_i^+$;\; $b_i = b_i^+$ for all $i\in [n]$;\; $S = S^+$
            \hfill $\diamond$ Transfer buffered gradients to main table
            \label{line:transfer_table}
    \ENDWHILE
\end{algorithmic}
\end{algorithm}
\begin{algorithm}[H]\caption{Worker $i$'s subroutine}
    \label{algo:Ringleader_worker}
    \begin{algorithmic}[1]
        \STATE \textbf{Input:} Model $x$
        \WHILE{True}
            \STATE Compute $g_i = \nabla f_i(x; \xi_i)$ using a freshly sampled data point $\xi_i \sim \cD_i$ 
            \STATE Send $g_i$ to the server 
        \ENDWHILE
    \end{algorithmic}
\end{algorithm}
\end{minipage}
\end{figure}
\subsection{Detailed Description}
\paragraph{Initialization.}
The algorithm begins by broadcasting the initial point $x^0 \in \R^d$ to all workers, which then start executing the worker subroutine (\Cref{algo:Ringleader_worker}).
Each worker continuously computes stochastic gradients at its current point and sends them to the server until instructed to stop, at which point the server provides a new point to resume from.
This design ensures that workers remain fully utilized and never idle.

The server maintains a gradient table with entries $\{(G_i, b_i)\}_{i=1}^n$, all initialized to zero.
Here, $G_i$ accumulates gradients, while $b_i$ tracks the number of stochastic gradients received from worker $i$ to form proper minibatch estimators, with $b_i = 0$ for all $i \in [n]$ at the start.

Before Phase~1 begins, we also initialize the set $S = \emptyset$, which tracks which table entries contain at least one stochastic gradient.
Since no gradients have yet arrived, $S$ is empty.
\paragraph{Phase~1 --- Gradient Collection.}
In this phase, the server continuously receives stochastic gradients from the workers and stores them in the gradient table $\{(G_i, b_i)\}_{i=1}^n$.
We denote by $g_j^k$ the stochastic gradient sent by worker $j$ at iteration $k$, which is computed at point $x^{k-\delta_j^k}$ using an i.i.d. sample $\xi_j \sim D_j$.

We do not specify how the delays $\delta_j^k$ evolve, since this information is not needed to run the algorithm: whenever necessary, a delay can be obtained as the difference between the current model iterate and the iterate at which the stochastic gradient was computed.
The delays will only play a role in the analysis, not in the execution of the method.

Upon receiving $g_j^k$ from worker $j$ (Line~\ref{line:receive_ph1}), the server updates the corresponding table entry and the stochastic gradient counter as follows (Line~\ref{line:update_ph1})
$$
    G_j = G_j + g_j^k ~, \quad b_j = b_j + 1 ~, \quad S = S \cup \{j\} ~.
$$
This process continues until $S = [n]$, i.e., the server has collected at least one gradient from every worker.
No model updates are performed during this phase, and workers do not receive new points; hence, all stochastic gradients from a given worker are computed at the same point.
\paragraph{Phase~2 --- Sequential Updates.}
In this phase, the server performs exactly one model update for each worker $i$, for a total of $n$ updates.
Phase~2 starts with the last worker that completed Phase~1, i.e., the worker whose gradient made the table complete.
The server first computes an update by averaging the accumulated stochastic gradients in the table $\{(G_i, b_i)\}_{i=1}^n$ and taking a descent step with this estimate (Line~\ref{line:step_last_worker}).
The updated model is then sent to this worker (Line~\ref{line:broadcast_last_worker}), the worker is removed from the set $S$, and the iteration counter is incremented (Line~\ref{line:counter_last_worker}).

Next, the server must update the remaining $n-1$ workers.
These updates are performed sequentially as soon as each worker finishes its current computation.
During this waiting period, new gradients may arrive from workers not in $S$---e.g. for example, the last updated worker may send another stochastic gradient before the other workers complete their computation.
Since discarding these gradients would waste information, they are instead buffered for the next round.
\paragraph{Temporary Table Management.}
To achieve this, the server maintains a temporary table $\{(G_i^+, b_i^+)\}_{i=1}^n$, initialized to zero (Line~\ref{line:buffer_initialization}), together with a set $S^+$ that records which workers have contributed to the table.
Whenever a gradient arrives from a worker not in $S$, it is stored in the temporary table (Line~\ref{line:buffer}).

If instead the gradient comes from a worker $j \in S$---i.e., one of the workers whose model we still need to update---the server first updates the main table $\{(G_i, b_i)\}_{i=1}^n$ with this new stochastic gradient (Line~\ref{line:update_table_ph2}).
It then performs a model update by again averaging the accumulated stochastic gradients in the table (Line~\ref{line:step_ph2}), broadcasts the new model to worker $j$ (Line~\ref{line:broadcast_ph2}), increments the iteration counter, and removes $j$ from the set $S = S \setminus \{j\}$ (Line~\ref{line:counter_ph2}).
\paragraph{Preparing for the Next Round.}
Once all workers in $S$ have been updated and Phase~2 is complete ($S=\emptyset$), the server prepares for the next round by copying the contents of the temporary table $\{(G_i^+, b_i^+)\}_{i=1}^n$ into the main table $\{(G_i, b_i)\}_{i=1}^n$ (Line~\ref{line:transfer_table}).
The set $S$ is also reset to $S = S^+$, since these workers already contributed gradients at their updated models.
Entries in the main table corresponding to workers not in $S^+$ remain zero, as the temporary table was initialized with zeros at the start of Phase~2 (Line~\ref{line:buffer_initialization}).

The server can now proceed to the next round by returning to Phase~1 and beginning a new gradient collection phase.
\subsection{Delay Control Analysis}
The structure of \algn, with its two phases, is specifically designed to prevent the unbounded delays that arise in standard asynchronous methods.
To understand why this works, consider that in each round we perform exactly $n$ updates---one per worker---before moving to the next round.
This ensures that no worker can fall more than one full round behind the current iteration.
The precise bound on delays is given in the following lemma
%
\begin{boxedlemma}[Bounded Delay]\label{lemma:delay}
   In \algn (\Cref{algo:Ringleader}), the delays $\delta_i^k$ satisfy
   $$
       \delta_i^k \le 2n-2~,
   $$
   for any worker $i \in [n]$ and any iteration $k \geq 0$.
\end{boxedlemma}
\begin{proof}
    We prove this by analyzing the structure of \algn.
    The algorithm operates in rounds, where each round consists of Phase~1 (gradient collection) followed by Phase~2 (sequential updates). In each Phase~2, the server performs exactly $n$ updates, one for each worker.
    Phase~2 begins at iterations $0, n, 2n, 3n, \ldots$, i.e., at multiples of $n$.
    \paragraph{First round (iterations $0$ to $n-1$):}
    Initially, all workers compute gradients at the point $x^0$, so during iterations $0, 1, \ldots, n-1$, the server receives gradients computed at $x^0$.
    For any iteration $k$ in this range, the server processes stochastic gradients computed at point $x^{k-\delta_i^k}$, so $\delta_i^k = k \le n-1$.
    Thus, delays simply increment during this first Phase~2.

    At the end of this round, each worker $i$ has received a new model $x^j$ for some $j \in \{1, 2, \ldots, n\}$, and these update iterations are distinct across workers.
    \paragraph{Second round (iterations $n$ to $2n-1$):}
    At the start of the second Phase~2 (at iteration $n$), the gradient table contains gradients computed at points $x^{n-\delta_i^n}$ for worker $i$, where $\delta_i^n \in \{0, 1, \ldots, n-1\}$.
    These delay values are distinct across workers since each worker received its update at a different iteration in the previous round.
    
    During iterations $n$ to $2n-1$, these delays increase by $1$ at each iteration for the same reason as in the first Phase~2, giving $\delta_i^{2n-1} \in \{n-1, n, \ldots, 2n-2\}$ at the end of this round.
    At the same time, all workers receive new points to compute gradients from, so during the next Phase~2, the delays will again be distinct for all workers and in $\{0, 1, \ldots, n-1\}$.
    \paragraph{General pattern:}
    By induction, at the beginning of each round starting at iteration $cn$ (for integer $c \ge 1$), the delays $\delta_i^{cn}$ take distinct values in $\{0, 1, \ldots, n-1\}$.
    During each Phase~2, these delays increase by at most $n-1$, giving the bound
    $$
        \delta_i^k \le (n-1) + (n-1) = 2n-2~.
    $$
\end{proof}
\subsection{Comparison to \iasgdtitle}
Our method is a delay-controlled version of \iasgd.
We can recover \iasgd by removing Phase~1 (gradient collection) and Phase~2 (structured updates), and thus perform updates naively---immediately upon gradient arrival.
In contrast, our algorithm operates in structured rounds, performing exactly one update per worker in each round, which provides the crucial delay control that \iasgd lacks.

In \iasgd, delays for slow workers can grow unboundedly because the server continuously updates the model using gradients from fast workers, causing slow workers to fall increasingly behind.
Our method prevents this issue by buffering the gradients from fast workers rather than immediately applying these gradients, to ensure that all workers receive updated models within $n$ subsequent iterations.
\subsection{Comparison to \maleniatitle}
\malenia also operates as an algorithm with two phases.
In Phase~1, \malenia collects gradients using a similar method to our approach, but uses a different termination condition \eqref{eq:malenia_condition} that requires knowledge of the noise parameter $\sigma$ and the target stationarity level $\varepsilon$, making it impractical.
In Phase~2, \malenia performs a \textit{synchronous} update by averaging all collected gradients and broadcasting the new model to \textit{all} workers simultaneously before returning to Phase~1.
This synchronization forces \malenia to discard ongoing gradient computations from workers that are active during the broadcast.

In contrast, our method performs Phase~2 \textit{asynchronously}: we update workers sequentially as they complete their current gradient computations, which ensures that no computational work is wasted.

Regarding \malenia's termination condition \eqref{eq:malenia_condition}, in \Cref{sec:malenia_param_free} we demonstrate that this condition can be replaced with our simpler requirement of obtaining at least one gradient from every worker.
With this modification, \malenia remains optimal in the fixed-time regime \eqref{eq:fixed_time} while becoming parameter-free, which eliminates the need for prior knowledge of $\sigma$ and $\varepsilon$.
Under this parameter-free variant, the only difference between \malenia and \algn lies in Phase~2: we perform updates asynchronously without discarding gradients, while \malenia operates synchronously.
\section{Theoretical Results}
Before presenting the theoretical results, we first write the algorithm in a compact form.
The gradients for each worker in the table are all computed at the same point; for worker $i$ at iteration $k$, the point is $x^{k-\delta_i^k}$.
The update rule can be written compactly as
$$
    x^{k+1} = x^k - \gamma \bar g^k~,
$$
where the gradient estimator $\bar g^k$ is defined by
$$
    \bar g^k
    \eqdef \frac{1}{n} \sum_{i=1}^n \bar g_i^k
    \eqdef \frac{1}{n} \sum_{i=1}^n \frac{1}{b_i^k} \sum_{j=1}^{b_i^k} g_i^{k,j} ~.
$$
Since multiple gradients may be received from the same worker, we denote by $g_i^{k,j}$ the $j$-th gradient from worker $i$ at iteration $k$.
Here the index $j$ corresponds to the i.i.d. sampled data point, and more concretely
$$
    g_i^{k,j} \eqdef \nabla f_i\(x^{k-\delta_i^k};\xi_i^{k-\delta_i^k, j}\)~.
$$
The quantity $b_i^k$ denotes the number of gradients from worker $i$ stored in the table at iteration $k$, i.e., the value of $b_i$ in Lines~\ref{line:step_last_worker} and \ref{line:step_ph2}.
Thus, the pair $(b_i^k, \delta_i^k)$ fully determines the method's behavior at iteration $k$.

Note that the sequence $\{b_i^k\}$ depends only on the computation times $\{\tau_i\}$ and the algorithm design (i.e., the stopping rule for collecting gradients).
Once these are fixed, all $b_i^k$ for every $i \in [n]$ and iteration $k$ are determined.
Crucially, the values of $b_i^k$ do not depend on the method's hyperparameters $\gamma$, $x^0$, or on the variance parameter $\sigma$ or the stationarity level $\varepsilon$.
\subsection{Iteration Complexity}
Our convergence analysis follows the structured approach employed by \citet{maranjyan2025ringmaster}, which decomposes the proof into two key components: a descent lemma that tracks the progress of the objective function and a residual estimation lemma that controls the accumulated delays in the system.

We begin by establishing notation for the harmonic means of the batch sizes across rounds:
$$
   B^k = \( \frac{1}{n} \sum_{i=1}^n \frac{1}{b_i^k} \)^{-1}, \quad \text{and} \quad B = \inf_{k \ge 0} B^k.
$$
Note that $B \geq 1$, since by the algorithm's design each $b_i^k \geq 1$.
A sharper bound on $B$ will be established later in \Cref{lemma:time_for_n_iter}.

The first lemma quantifies how much the objective function decreases at each iteration, accounting for both the standard gradient descent progress and the additional complexities introduced by asynchronous updates.
\begin{boxedlemma}[Descent Lemma; Proof in \Cref{proof:descent}]\label{lemma:descent}
    Under Assumptions~\ref{ass:stochastic_variance_bounded} and~\ref{ass:lipschitz_constant}, if the stepsize in \Cref{algo:Ringleader} satisfies $\gamma \le \nicefrac{1}{4L}$, then the following inequality holds
    \begin{align*}
        \E{f\(x^{k+1}\)}
        &\le \E{f\(x^{k}\)}
            - \frac{\gamma}{2} \E{ \sqnorm{\nabla f\(x^{k}\)} }
            - \frac{\gamma}{4} \E{\sqnorm{\frac{1}{n} \sum_{i=1}^n \nabla f_i\(x^{k-\delta_i^k}\)}} \\
            &\quad + \frac{\gamma L^2}{2n} \sum_{i=1}^n \E{\sqnorm{x^{k} - x^{k-\delta_{i}^k}}}
                    + \frac{3\gamma^2 L \sigma^2}{2B} \\
            &\quad + \gamma^2 L \sum_{\ell = k-(k \bmod n)}^{k-1} \E{\sqnorm{\frac{1}{n}\sum_{i=1}^n \nabla f_i\(x^{\ell-\delta_i^\ell}\)}}.
    \end{align*}
\end{boxedlemma}
This descent lemma shares a similar structure with its counterpart in the homogeneous setting analyzed by \citet{maranjyan2025ringmaster}, but with a crucial additional term.
The final summation term in the upper bound captures the effect of using stale gradients from the gradient table---a phenomenon we refer to as ``table delay".
This term is absent in the homogeneous case because no gradient table is maintained.
Indeed, when $n=1$, our setting reduces to the homogeneous case, the gradient table becomes unnecessary, and this additional term vanishes, recovering the original descent lemma established by \citet{maranjyan2025ringmaster}.

Next, similar to the work by \citet{maranjyan2025ringmaster}, we derive a lemma to bound the term involving the difference between current and old points.
\begin{boxedlemma}[Residual Estimation; Proof in \Cref{proof:residual}]
    \label{lemma:residual}
    Under \Cref{ass:stochastic_variance_bounded}, the iterates of \algn (\Cref{algo:Ringleader}) with stepsize $\gamma \le \nicefrac{1}{4nL}$ satisfy the following bound
    \begin{equation*}
    \frac{1}{K} \sum_{k=0}^{K-1} \frac{1}{n} \sum_{i=1}^n \E{\sqnorm{x^{k} - x^{k-\delta_i^k}}}
        \le \frac{2\gamma n}{LK}\sum_{k=0}^{K-1} \E{\sqnorm{ \frac{1}{n}\sum_{j=1}^n \nabla f_j\(x^{k-\delta_j^k}\) }}
            + \frac{2\gamma \sigma^2}{LB}~.
    \end{equation*}
\end{boxedlemma}
Finally, we get the iteration complexity combining these two lemmas.
\begin{boxedtheorem}[Iteration Complexity]\label{theorem:convergence}
    Under Assumptions \ref{ass:stochastic_variance_bounded}, \ref{ass:lipschitz_constant}, and \ref{ass:lower_bound}, let the stepsize in \algn (\Cref{algo:Ringleader}) be
    $$
        \gamma = \min \left\{\frac{1}{8nL}, \frac{\varepsilon B}{10 L \sigma^2 } \right\}.
    $$
    Then,
    $$
        \frac{1}{K}\sum_{k=0}^{K-1} \E{ \sqnorm{\nabla f\(x^{k}\)} } 
        \le \varepsilon~,
    $$
    for
    $$
        K 
        \ge \frac{32 nL\Delta }{\varepsilon} + \frac{40 L\Delta \sigma^2}{B\varepsilon^2}
        = \cO\left(\frac{nL\Delta}{\varepsilon} \left(1 + \frac{\sigma^2}{Bn\varepsilon}\right)\right).
    $$
\end{boxedtheorem}
\begin{proof}    
    We start by averaging the inequality from \Cref{lemma:descent} over $K$ iterations and dividing both sides by $\nicefrac{\gamma}{2}$
    \begin{align*}
        \frac{1}{K}\sum_{k=0}^{K-1} & \E{ \sqnorm{\nabla f\(x^{k}\)} } 
            + \frac{1}{2K}\sum_{k=0}^{K-1} \E{\sqnorm{\frac{1}{n} \sum_{i=1}^n \nabla f_i\(x^{k-\delta_i^k}\)}} \\
        &\le \frac{2\Delta}{\gamma K}
            + \frac{3 \gamma L \sigma^2}{B} \\
            &\quad + \frac{L^2}{n} \frac{1}{K}\sum_{k=0}^{K-1} \sum_{i=1}^n \E{\sqnorm{x^{k} - x^{k-\delta_{i}^k}}} \\
            &\quad + 2 \gamma L \frac{1}{K}\sum_{k=0}^{K-1} \sum_{\ell = k- (k \bmod n)}^{k-1} \E{\sqnorm{\frac{1}{n}\sum_{i=1}^n \nabla f_i\(x^{\ell-\delta_i^\ell}\)}}.
    \end{align*}
    We now bound the third term on the right using \Cref{lemma:residual}
    \begin{align*}
        \frac{1}{K}\sum_{k=0}^{K-1} & \E{ \sqnorm{\nabla f\(x^{k}\)} } 
            + \frac{1}{2K}\sum_{k=0}^{K-1} \E{\sqnorm{\frac{1}{n} \sum_{i=1}^n \nabla f_i\(x^{k-\delta_i^k}\)}} \\
        &\le \frac{2\Delta}{\gamma K}
            + \frac{3 \gamma L \sigma^2}{B}
            + \frac{2\gamma L \sigma^2}{B} \\
            &\quad + 2\gamma L n \frac{1}{K}\sum_{k=0}^{K-1} \E{\sqnorm{ \frac{1}{n}\sum_{j=1}^n \nabla f_j\(x^{k-\delta_j^k}\) }} \\
            &\quad + 2 \gamma L \frac{1}{K}\sum_{k=0}^{K-1} \sum_{\ell = k- (k \bmod n)}^{k-1} \E{\sqnorm{\frac{1}{n}\sum_{i=1}^n \nabla f_i\(x^{\ell-\delta_i^\ell}\)}} \\
        &\le \frac{2\Delta}{\gamma K}
            + \frac{5 \gamma L \sigma^2}{B} \\
            &\quad + 2\gamma L n\frac{1}{K}\sum_{k=0}^{K-1} \E{\sqnorm{ \frac{1}{n}\sum_{j=1}^n \nabla f_j\(x^{k-\delta_j^k}\) }} \\
            &\quad + 2 \gamma Ln \frac{1}{K}\sum_{k=0}^{K-1} \E{\sqnorm{\frac{1}{n}\sum_{i=1}^n \nabla f_i\(x^{k-\delta_i^k}\)}} .
    \end{align*}
    Now, using the bound $\gamma \le \nicefrac{1}{8nL}$, we obtain
    \begin{equation*}
        \frac{1}{K}\sum_{k=0}^{K-1} \E{ \sqnorm{\nabla f\(x^{k}\)} }
        \le \frac{2\Delta}{\gamma K}
            + \frac{5 \gamma L \sigma^2}{B}~.
    \end{equation*}
    Finally, plugging in the stepsize and the expression for $K$ ensures the right-hand side is bounded by $\varepsilon$.
\end{proof}
For parallel and asynchronous methods, iteration complexity is less important than time complexity.
What truly matters is how quickly we can finish training.
We are willing to perform more iterations and extra computation if it means completing the process faster.
Having established the iteration complexity, we now turn to the time complexity.
\subsection{Time Complexity}
Since the algorithm operates in rounds with $n$ steps per round, and its iteration complexity is already known, it remains to determine the duration of each round.
We have the following lemma
\begin{boxedlemma}\label{lemma:time_for_n_iter}
    Each block of $n$ consecutive iterations (each round) of \Cref{algo:Ringleader} takes at most $2\tau_n$ seconds.
    Moreover, we have
    $$
        B \ge \frac{\tau_n}{2} \left( \frac{1}{n} \sum_{i=1}^n \tau_i \right)^{-1} = \frac{\tau_n}{2\tauavg}~.
    $$
\end{boxedlemma}
\begin{proof}
   The upper bound of $2\tau_n$ follows from the structure of the algorithm, which consists of two phases.
   In the first phase, the server waits until all workers complete at least one gradient computation, which takes at most $\tau_n$ seconds.
   In the second phase, the server applies the received gradients and waits for all ongoing computations to finish---which again takes at most $\tau_n$ seconds.
   Thus, the total time for $n$ iterations is bounded by $2\tau_n$.

   We now prove the second part of the lemma.
   Recall that
   $$
      B = \inf_{k \ge 0} B^k = \inf_{k \ge 0} \left( \frac{1}{n} \sum_{i=1}^n \frac{1}{b_i^k} \right)^{-1} \ge \left( \frac{1}{n} \sum_{i=1}^n \frac{1}{b_i} \right)^{-1},
   $$
   where we define
   $$
      b_i = \inf_{k \ge 0} b_i^k~.
   $$
   We are interested in the number of gradients stored in the table at iteration $k$.
   This count includes gradients computed during Phase~1 plus one additional gradient from Phase~2 (except for the worker that finished Phase~1 last).

   Since every worker needs to compute at least one gradient during Phase~1, the slowest worker will take $\tau_n$ seconds to complete single gradient computation.
   During this $\tau_n$-second interval, faster workers $i < n$ may still be finishing gradients from the previous round's Phase~2 before starting their Phase~1 computations for the current round.

   After completing any remaining Phase~2 work (which takes at most $\tau_i$ seconds), worker $i$ has at least $\tau_n - \tau_i$ seconds remaining to compute additional gradients for the current round's Phase~1.
   The number of gradients that worker $i$ can compute satisfies
   $$
      b_i \ge \max\left\{1, \left\lceil \frac{\tau_n - \tau_i}{\tau_i} \right\rceil \right\} \ge \max\left\{1, \frac{\tau_n}{\tau_i} - 1 \right\}.
   $$
   For workers $i$ where $\tau_n \geq 2\tau_i$, we have
   $$
      \frac{\tau_n}{\tau_i} - 1 \ge \frac{\tau_n}{2\tau_i}~,
   $$
   and hence
   $$
      b_i \ge \frac{\tau_n}{2\tau_i}~.
   $$
   Plugging this bound into the expression for $B$ gives the claimed result.
\end{proof}
Based on this lemma, we derive the time complexity guarantee of our algorithm
\begin{boxedtheorem}\label{thm:time_complexity}
    Let Assumptions \ref{ass:lipschitz_constant}, \ref{ass:lower_bound}, and \ref{ass:stochastic_variance_bounded} hold.
    Let the stepsize in \algn (\Cref{algo:Ringleader}) be
    $\gamma = \min \left\{\nicefrac{1}{8nL}, \nicefrac{\varepsilon B}{10 L \sigma^2 } \right\}$.
    Then, under the \emph{fixed computation model} \eqref{eq:fixed_time}, \algn achieves the optimal time complexity
    $$
        \cO\(\frac{L\Delta}{\varepsilon}\( \tau_n + \tauavg\frac{\sigma^2}{n\varepsilon} \) \) .
    $$
\end{boxedtheorem}
\begin{proof}
    We start with the iteration complexity from \Cref{theorem:convergence}
    $$
        K 
        \ge \frac{32 nL\Delta }{\varepsilon} + \frac{40 L\Delta \sigma^2}{B\varepsilon^2}
        = \cO\(\frac{nL\Delta}{\varepsilon} \(1 + \frac{\sigma^2}{B n\varepsilon}\) \).
    $$
    The time to do $n$ steps is at most $2\tau_n$ form \Cref{lemma:time_for_n_iter}.
    Then the time complexity is 
    $$
        2\tau_n \times \frac{K}{n} 
        = \cO\(\tau_n\frac{L\Delta}{\varepsilon} \( 1 + \frac{\sigma^2}{Bn\varepsilon} \) \) .
    $$
    It remains to put $B \ge \nicefrac{\tau_n}{2\tauavg}$ from \Cref{lemma:time_for_n_iter}.
\end{proof}
The obtained time complexity consists of two key terms that illuminate the algorithm's behavior.
The first term depends on the slowest device, which is fundamental since all devices must contribute to solving the problem.
The second term, however, involves $\tauavg$ rather than $\tau_n$ as in \naiveminibatch (see \Cref{table:fixedtime})---this substitution captures the core benefit of asynchronous execution.
Specifically, this advantage becomes pronounced when $\sigma$ is relatively large.
Intuitively, in high-noise regimes, collecting many gradients from workers is essential for convergence, and asynchronous methods can leverage faster workers more effectively.
Conversely, in low-noise settings, fewer gradient evaluations suffice for good performance, making \naiveminibatch already quite effective and rendering the additional complexity of asynchrony unnecessary.
\begin{remark}
    The optimality claim for \Cref{thm:time_complexity} holds when the smoothness-type constant~$L$ in \Cref{ass:lipschitz_constant} is within a constant factor of the smoothness~$L_f$ used to derive the lower bound \citep{tyurin2024optimal} (\Cref{table:fixedtime}).
\end{remark}
Under this condition, the resulting time complexity matches the lower bound of \citet{tyurin2024optimal}, making \algn the first asynchronous algorithm to achieve optimality under heterogeneous data.
\section{Experiments}
\label{sec:experiments}
To validate our theoretical results we perform a toy simulation.

We consider image classification on MNIST \citep{lecun2010mnist} and on Fashion-MNIST \citep{xiao2017fashion} with standard normalization for both datasets.
To enable equal client sizes, we first trim the dataset so that the total number of examples is divisible by the number of clients $n=100$.
To obtain heterogeneous datasets across clients, we then partition the trimmed set using an \emph{equal-size Dirichlet} procedure with concentration parameter $\alpha=0.1$ \citep{yurochkin2019bayesian}.
For each client $j\in[n]$, we draw proportions $p_j \sim \mathrm{Dirichlet}(\alpha,\ldots,\alpha)$ over the classes and compute a rounded class-allocation vector whose entries sum exactly to $\nicefrac{N}{n}$, where $N$ is the trimmed dataset size.
This creates non-IID data where each client has a skewed distribution over classes (with $\alpha=0.1$, clients typically observe only 1-2 classes frequently).

When assigning samples, we take the requested number from each class pool for client $j$.
If a class pool does not have enough remaining examples to meet the requested amount, the client receives whatever is left from that class and the shortfall is \emph{topped up} using samples from other classes that still have available examples.

Our model is a two-layer MLP $\mathrm{Linear}(d,128)\!\to\!\mathrm{ReLU}\!\to\!\mathrm{Linear}(128,10)$ trained with mean cross-entropy.
Stochastic gradients at the clients use a minibatch of size $4$, while reported gradient norms are computed on the \emph{full} dataset.

We emulate heterogeneous compute by assigning each client $i$ a base delay and a random jitter:
$$
    \tau_i \;=\; i \;+\; |\eta_i|~, \qquad 
    \eta_i \sim \cN(0,i)~, 
    \quad \text{for all } i\in [n]~.
$$  
For each method we tune the stepsize $\gamma$ within a fixed wall-clock budget.
We sweep  
$$
    \gamma \in \{0.001,\,0.005,\,0.01,\,0.02,\,0.05,\,0.1,\,0.2,\,0.5,\,1.0,\,2.0\}~,
$$  
and then fix the best $\gamma$ per method for evaluation.

We report the full-batch gradient-norm squared $\|\nabla f(x^k)\|^2$, versus wall-clock time.
Each method is run $30$ times over different seeds.
We report the \emph{median} with \emph{interquartile range (IQR)}.
To reduce high-frequency noise, we apply a centered moving-average smoothing to the aggregated curves (post-aggregation), while keeping the initial point unchanged.

\begin{figure}[t]
    \centering
    \begin{minipage}{0.48\linewidth}
        \centering
        \includegraphics[width=\linewidth]{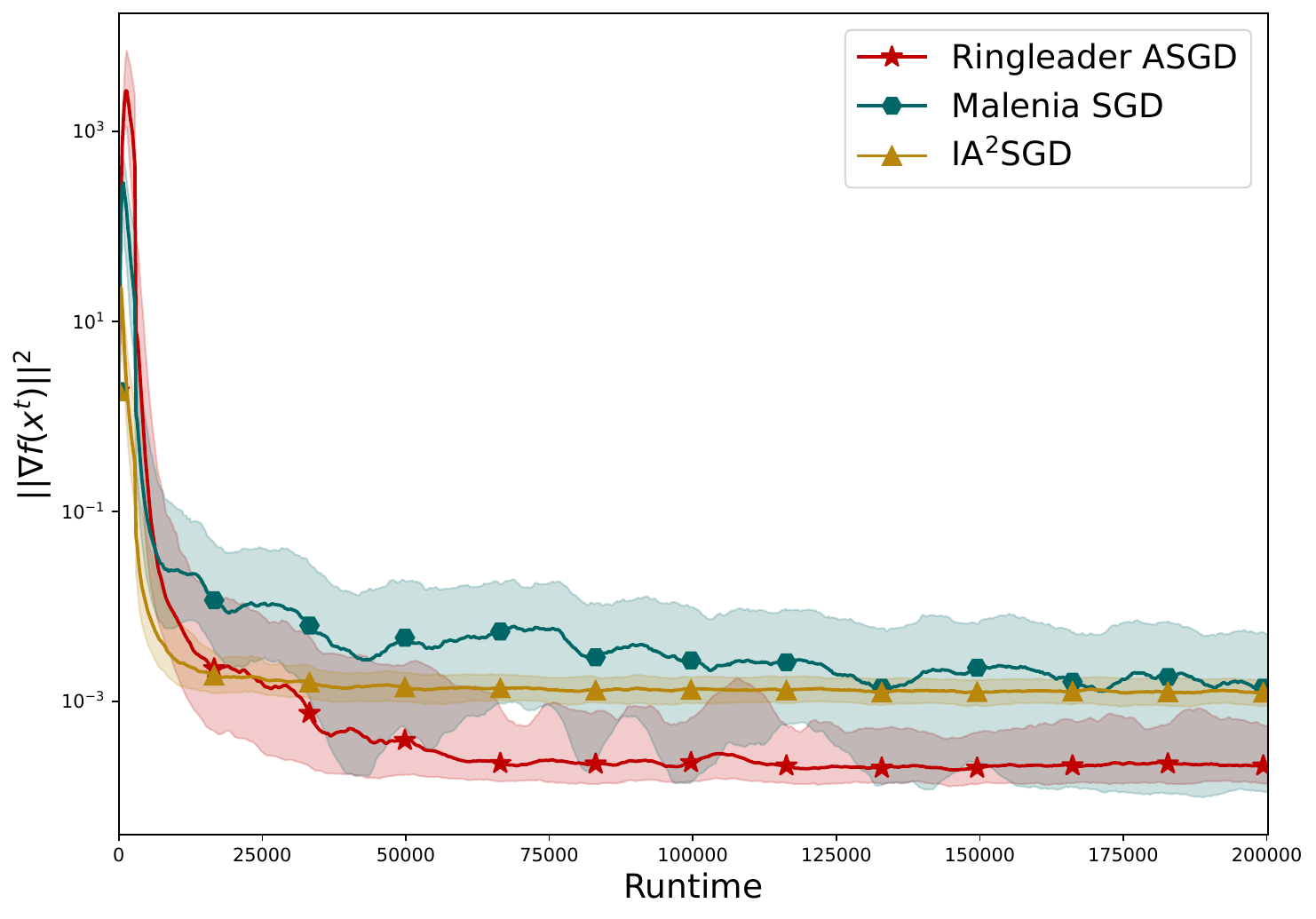}
        \caption*{(a) MNIST}
    \end{minipage}
    \hfill
    \begin{minipage}{0.48\linewidth}
        \centering
        \includegraphics[width=\linewidth]{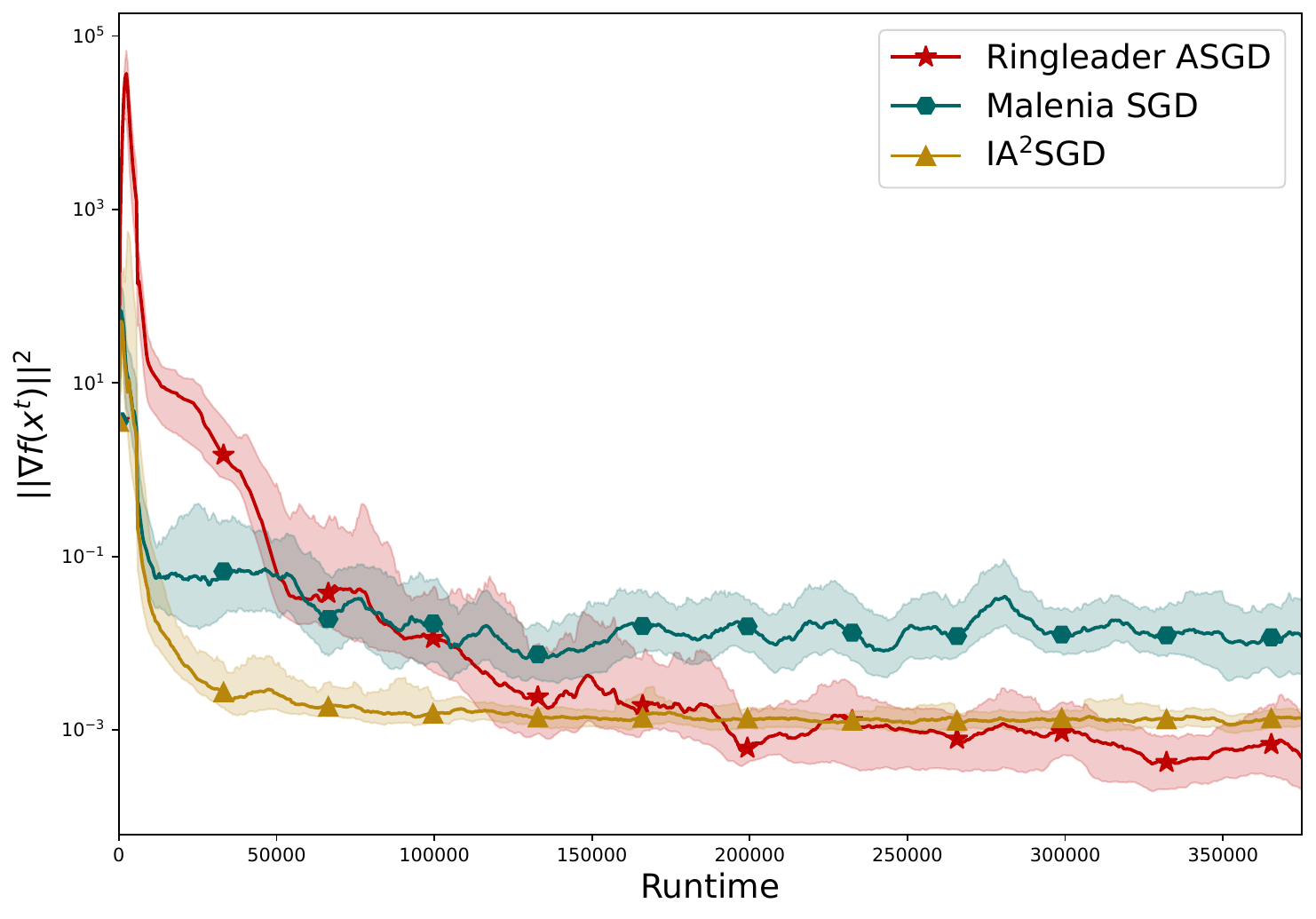}
        \caption*{(b) Fashion-MNIST}
    \end{minipage}
    \caption{
        Convergence comparison showing median gradient norm squared $\|\nabla f(x^k)\|^2$ (solid lines) with interquartile ranges (shaded regions) versus wall-clock time, averaged over 30 random seeds.
        \textbf{Setup:} Two-layer MLP with architecture $\mathrm{Linear}(d,128) \to \mathrm{ReLU} \to \mathrm{Linear}(128,10)$ trained on (a) MNIST and (b) Fashion-MNIST datasets.
        \textbf{Client delays:} Heterogeneous delays simulated as $\tau_i = i + |\eta_i|$ where $\eta_i \sim \mathcal{N}(0,i)$ for client $i \in [n]$, where we choose $n=100$.
        \textbf{Results:} With optimally tuned stepsizes, \algn achieves faster convergence than both \malenia and \iasgd, despite \algn and \malenia having equivalent time complexity guarantees.
    }
    \label{fig:median}
\end{figure}

\Cref{fig:median} shows the results.
We observe that \algn converges faster compared to \malenia and \iasgd.  
Although theory suggests that \algn and \malenia have the same time complexity, in practice \algn benefits from the $n$ updates performed in Phase 2 instead of one synchronous update.
This design enables more optimization steps within the same wall-clock budget, which is especially advantageous when updates are sparse.
\section{Conclusion}
We have introduced \algn, the first asynchronous stochastic gradient method to achieve optimal time complexity under arbitrary data heterogeneity and arbitrarily heterogeneous computation times in distributed learning, without requiring similarity assumptions between workers' datasets.

Its core innovation is a two-phase structure within each round: the model is updated once per worker (for a total of $n$ updates), while a buffering mechanism manages gradient delays and preserves the efficiency of asynchronous execution.
By maintaining a gradient table and alternating between gradient collection and sequential updates, \algn prevents the unbounded delays common in naive asynchronous methods.
Every gradient received by the server is either used in the current round or stored for future use, ensuring no computation is wasted.

Our analysis shows that \algn matches the optimal time complexity bounds established by \citet{tyurin2024optimal}.
In contrast to the optimal but synchronous \malenia method, \algn is asynchronous and requires no prior knowledge of problem parameters in the algorithm design, making it practical for real-world deployments.

Finally, with a minor modification, \algn also achieves optimality in the more general setting of arbitrarily varying computation times (\Cref{sec:arbitrary_time}).

\subsubsection*{Acknowledgments}
The research reported in this publication was supported by funding from King Abdullah University of Science and Technology (KAUST): i) KAUST Baseline Research Scheme, ii) CRG Grant ORFS-CRG12-2024-6460, and iii) Center of Excellence for Generative AI, under award number 5940.

\bibliography{bib}

@inproceedings{maranjyan2025ringmaster,
    title={{R}ingmaster {ASGD}: The First {A}synchronous {SGD} with Optimal Time Complexity}, 
    author={Artavazd Maranjyan and Alexander Tyurin and Peter Richtárik},
    year={2025},
    booktitle={International Conference on Machine Learning},
}

@inproceedings{maranjyan2025ata,
    title={{ATA}: Adaptive Task Allocation for Efficient Resource Management in Distributed Machine Learning}, 
    author={Maranjyan, Artavazd and Saad, El Mehdi and Richt{\'a}rik, Peter and Orabona, Francesco},
    year={2025},
    booktitle={International Conference on Machine Learning},
}

@inproceedings{maranjyan2025mindflayer,
  title={MindFlayer {SGD}: Efficient Parallel {SGD} in the Presence of Heterogeneous and Random Worker Compute Times},
  author={Artavazd Maranjyan and Omar Shaikh Omar and Peter Richt{\'a}rik},
  booktitle={The 41st Conference on Uncertainty in Artificial Intelligence},
  year={2025},
  url={https://openreview.net/forum?id=RNpvu3MSvm}
}

@article{maranjyan2025gradskip,
  title={{GradSkip}: Communication-Accelerated Local Gradient Methods with Better Computational Complexity},
  author={Artavazd Maranjyan and Mher Safaryan and Peter Richt{\'a}rik},
  journal={Transactions on Machine Learning Research},
  issn={2835-8856},
  year={2025},
  url={https://openreview.net/forum?id=6R3fRqFfhn},
}

@inproceedings{zheng2017asynchronous,
  title={Asynchronous stochastic gradient descent with delay compensation},
  author={Zheng, Shuxin and Meng, Qi and Wang, Taifeng and Chen, Wei and Yu, Nenghai and Ma, Zhi-Ming and Liu, Tie-Yan},
  booktitle={International Conference on Machine Learning},
  pages={4120--4129},
  year={2017},
  organization={PMLR}
}

@article{xie2019asynchronous,
  title={Asynchronous federated optimization},
  author={Xie, Cong and Koyejo, Sanmi and Gupta, Indranil},
  journal={arXiv preprint arXiv:1903.03934},
  year={2019}
}

@article{fraboni2023general,
  title={A general theory for federated optimization with asynchronous and heterogeneous clients updates},
  author={Fraboni, Yann and Vidal, Richard and Kameni, Laetitia and Lorenzi, Marco},
  journal={Journal of Machine Learning Research},
  volume={24},
  number={110},
  pages={1--43},
  year={2023}
}

@article{lian2015asynchronous,
  title={Asynchronous parallel stochastic gradient for nonconvex optimization},
  author={Lian, Xiangru and Huang, Yijun and Li, Yuncheng and Liu, Ji},
  journal={Advances in Neural Information Processing Systems},
  volume={28},
  year={2015}
}

@article{tyurin2024tighttimecomplexitiesparallel,
  title={Tight time complexities in parallel stochastic optimization with arbitrary computation dynamics},
  author={Tyurin, Alexander},
  journal={arXiv preprint arXiv:2408.04929},
  year={2024}
}

@article{wang2022asyncfeded,
  title={{AsyncFedED}: Asynchronous federated learning with Euclidean distance based adaptive weight aggregation},
  author={Wang, Qiyuan and Yang, Qianqian and He, Shibo and Shi, Zhiguo and Chen, Jiming},
  journal={arXiv preprint arXiv:2205.13797},
  year={2022}
}

@article{wang2022asynchronous,
  title={Asynchronous federated learning over wireless communication networks},
  author={Wang, Zhongyu and Zhang, Zhaoyang and Tian, Yuqing and Yang, Qianqian and Shan, Hangguan and Wang, Wei and Quek, Tony QS},
  journal={IEEE Transactions on Wireless Communications},
  volume={21},
  number={9},
  pages={6961--6978},
  year={2022},
  publisher={IEEE}
}

@article{tsitsiklis1986distributed,
  title={Distributed asynchronous deterministic and stochastic gradient optimization algorithms},
  author={Tsitsiklis, John and Bertsekas, Dimitri and Athans, Michael},
  journal={IEEE Transactions on Automatic Control},
  volume={31},
  number={9},
  pages={803--812},
  year={1986},
  publisher={IEEE}
}

@article{feyzmahdavian2023asynchronous,
  title={Asynchronous iterations in optimization: New sequence results and sharper algorithmic guarantees},
  author={Feyzmahdavian, Hamid Reza and Johansson, Mikael},
  journal={Journal of Machine Learning Research},
  volume={24},
  number={158},
  pages={1--75},
  year={2023}
}

@article{shoeybi2019megatron,
  title={Megatron-{LM}: Training Multi-Billion Parameter Language Models Using Model Parallelism},
  author={Shoeybi, Mohammad and Patwary, Mostofa and Puri, Raul and LeGresley, Patrick and Casper, Jared and Catanzaro, Bryan},
  journal={arXiv preprint arXiv:1909.08053},
  year={2019}
}

@inproceedings{narayanan2021efficient,
  title={Efficient large-scale language model training on {GPU} clusters using {Megatron-LM}},
  author={Narayanan, Deepak and Shoeybi, Mohammad and Casper, Jared and LeGresley, Patrick and Patwary, Mostofa and Korthikanti, Vijay and Vainbrand, Dmitri and Kashinkunti, Prethvi and Bernauer, Julie and Catanzaro, Bryan and others},
  booktitle={Proceedings of the International Conference for High Performance Computing, Networking, Storage and Analysis},
  pages={1--15},
  year={2021}
}

@article{grattafiori2024llama,
  title={The Llama 3 herd of models},
  author={Grattafiori, Aaron and Dubey, Abhimanyu and Jauhri, Abhinav and Pandey, Abhinav and Kadian, Abhishek and Al-Dahle, Ahmad and Letman, Aiesha and Mathur, Akhil and Schelten, Alan and Vaughan, Alex and others},
  journal={arXiv preprint arXiv:2407.21783},
  year={2024}
}

@inproceedings{dean2012large,
 author = {Dean, Jeffrey and Corrado, Greg and Monga, Rajat and Chen, Kai and Devin, Matthieu and Mao, Mark and Ranzato, Marc\textquotesingle aurelio and Senior, Andrew and Tucker, Paul and Yang, Ke and Le, Quoc and Ng, Andrew},
 booktitle = {Advances in Neural Information Processing Systems},
 editor = {F. Pereira and C.J. Burges and L. Bottou and K.Q. Weinberger},
 pages = {},
 publisher = {Curran Associates, Inc.},
 title = {Large Scale Distributed Deep Networks},
 url = {https://proceedings.neurips.cc/paper_files/paper/2012/file/6aca97005c68f1206823815f66102863-Paper.pdf},
 volume = {25},
 year = {2012}
}

@article{li2014communication,
  title={Communication efficient distributed machine learning with the parameter server},
  author={Li, Mu and Andersen, David G and Smola, Alexander and Yu, Kai},
  journal={Advances in Neural Information Processing Systems},
  volume={27},
  year={2014}
}

@article{zhao2018federated,
  title={Federated learning with non-iid data},
  author={Zhao, Yue and Li, Meng and Lai, Liangzhen and Suda, Naveen and Civin, Damon and Chandra, Vikas},
  journal={arXiv preprint arXiv:1806.00582},
  year={2018}
}

@article{tan2022towards,
  title={Towards personalized federated learning},
  author={Tan, Alysa Ziying and Yu, Han and Cui, Lizhen and Yang, Qiang},
  journal={IEEE Transactions on Neural Networks and Learning Systems},
  volume={34},
  number={12},
  pages={9587--9603},
  year={2022},
  publisher={IEEE}
}

@article{blatt2007convergent,
  title={A convergent incremental gradient method with a constant step size},
  author={Blatt, Doron and Hero, Alfred O and Gauchman, Hillel},
  journal={SIAM Journal on Optimization},
  volume={18},
  number={1},
  pages={29--51},
  year={2007},
  publisher={SIAM}
}

@article{gurbuzbalaban2017convergence,
  title={On the convergence rate of incremental aggregated gradient algorithms},
  author={Gurbuzbalaban, Mert and Ozdaglar, Asuman and Parrilo, Pablo A},
  journal={SIAM Journal on Optimization},
  volume={27},
  number={2},
  pages={1035--1048},
  year={2017},
  publisher={SIAM}
}

@article{roux2012stochastic,
  title={A stochastic gradient method with an exponential convergence rate for finite training sets},
  author={Roux, Nicolas and Schmidt, Mark and Bach, Francis},
  journal={Advances in Neural Information Processing Systems},
  volume={25},
  year={2012}
}

@article{schmidt2017minimizing,
  title={Minimizing finite sums with the stochastic average gradient},
  author={Schmidt, Mark and Le Roux, Nicolas and Bach, Francis},
  journal={Mathematical Programming},
  volume={162},
  pages={83--112},
  year={2017},
  publisher={Springer}
}

@inproceedings{zhang2023no,
  title={No one idles: Efficient heterogeneous federated learning with parallel edge and server computation},
  author={Zhang, Feilong and Liu, Xianming and Lin, Shiyi and Wu, Gang and Zhou, Xiong and Jiang, Junjun and Ji, Xiangyang},
  booktitle={International Conference on Machine Learning},
  pages={41399--41413},
  year={2023},
  organization={PMLR}
}

@article{wang2024achieving,
  title={Achieving linear speedup in asynchronous federated learning with heterogeneous clients},
  author={Wang, Xiaolu and Li, Zijian and Jin, Shi and Zhang, Jun},
  journal={IEEE Transactions on Mobile Computing},
  year={2024},
  publisher={IEEE}
}

@inproceedings{glasgow2022asynchronous,
  title={Asynchronous distributed optimization with stochastic delays},
  author={Glasgow, Margalit R and Wootters, Mary},
  booktitle={International Conference on Artificial Intelligence and Statistics},
  pages={9247--9279},
  year={2022},
  organization={PMLR}
}

@misc{wang2025incremental,
  title={Incremental Aggregated Asynchronous {SGD} for Arbitrarily Heterogeneous Data},
  author={Xiaolu Wang and Yuchang Sun and Hoi To Wai and Jun Zhang},
  year={2025},
  url={https://openreview.net/forum?id=m3x4kDbYAK}
}

@article{vanli2018global,
  title={Global convergence rate of proximal incremental aggregated gradient methods},
  author={Vanli, N Denizcan and Gurbuzbalaban, Mert and Ozdaglar, Asuman},
  journal={SIAM Journal on Optimization},
  volume={28},
  number={2},
  pages={1282--1300},
  year={2018},
  publisher={SIAM}
}

@article{mania2017perturbed,
  title={Perturbed iterate analysis for asynchronous stochastic optimization},
  author={Mania, Horia and Pan, Xinghao and Papailiopoulos, Dimitris and Recht, Benjamin and Ramchandran, Kannan and Jordan, Michael I},
  journal={SIAM Journal on Optimization},
  volume={27},
  number={4},
  pages={2202--2229},
  year={2017},
  publisher={SIAM}
}

@article{agarwal2011distributed,
  title={Distributed delayed stochastic optimization},
  author={Agarwal, Alekh and Duchi, John C},
  journal={Advances in Neural Information Processing Systems},
  volume={24},
  year={2011}
}

@article{alahyane2025optimizing,
  title={Optimizing Asynchronous Federated Learning: A Delicate Trade-Off Between Model-Parameter Staleness and Update Frequency},
  author={Alahyane, Abdelkrim and Comte, C{\'e}line and Jonckheere, Matthieu and Moulines, {\'E}ric},
  journal={arXiv preprint arXiv:2502.08206},
  year={2025}
}

@inproceedings{GPT3,
 author = {Brown, Tom and Mann, Benjamin and Ryder, Nick and Subbiah, Melanie and Kaplan, Jared D and Dhariwal, Prafulla and Neelakantan, Arvind and Shyam, Pranav and Sastry, Girish and Askell, Amanda and Agarwal, Sandhini and Herbert-Voss, Ariel and Krueger, Gretchen and Henighan, Tom and Child, Rewon and Ramesh, Aditya and Ziegler, Daniel and Wu, Jeffrey and Winter, Clemens and Hesse, Chris and Chen, Mark and Sigler, Eric and Litwin, Mateusz and Gray, Scott and Chess, Benjamin and Clark, Jack and Berner, Christopher and McCandlish, Sam and Radford, Alec and Sutskever, Ilya and Amodei, Dario},
 booktitle = {Advances in Neural Information Processing Systems},
 editor = {H. Larochelle and M. Ranzato and R. Hadsell and M. F. Balcan and H. Lin},
 pages = {1877--1901},
 publisher = {Curran Associates, Inc.},
 title = {Language Models are Few-Shot Learners},
 url = {https://proceedings.neurips.cc/paper/2020/file/1457c0d6bfcb4967418bfb8ac142f64a-Paper.pdf},
 volume = {33},
 year = {2020}
}

@inproceedings{mcmahan2017communication,
  title={Communication-efficient learning of deep networks from decentralized data},
  author={McMahan, Brendan and Moore, Eider and Ramage, Daniel and Hampson, Seth and y Arcas, Blaise Aguera},
  booktitle={Artificial Intelligence and Statistics},
  pages={1273--1282},
  year={2017},
  organization={PMLR}
}

@book{nesterov2018lectures,
	title={Lectures on Convex Optimization},
	author={Nesterov, Yurii},
	volume={137},
	year={2018},
	publisher={Springer}
}

@article{lecun2010mnist,
  title={MNIST handwritten digit database},
  author={LeCun, Yann and Cortes, Corinna and Burges, CJ},
  journal={ATT Labs [Online]. Available: http://yann.lecun.com/exdb/mnist},
  volume={2},
  year={2010}
}

@inproceedings{konevcny2016federated,
  title	= {Federated Learning: Strategies for Improving Communication Efficiency},
  author	= {Jakub Konečný and H. Brendan McMahan and Felix X. Yu and Peter Richt\'{a}rik and Ananda Theertha Suresh and Dave Bacon},
  year	= {2016},
  booktitle	= {NIPS Workshop on Private Multi-Party Machine Learning}
}

@article{li2020federated,
  title={Federated optimization in heterogeneous networks},
  author={Li, Tian and Sahu, Anit Kumar and Zaheer, Manzil and Sanjabi, Maziar and Talwalkar, Ameet and Smith, Virginia},
  journal={Proceedings of Machine Learning and Systems},
  volume={2},
  pages={429--450},
  year={2020}
}

@article{kairouz2021advances,
  title={Advances and open problems in federated learning},
  author={Kairouz, Peter and McMahan, H Brendan and Avent, Brendan and Bellet, Aur{\'e}lien and Bennis, Mehdi and Bhagoji, Arjun Nitin and Bonawitz, Kallista and Charles, Zachary and Cormode, Graham and Cummings, Rachel and others},
  journal={Foundations and Trends{\textregistered} in Machine Learning},
  volume={14},
  number={1--2},
  pages={1--210},
  year={2021},
  publisher={Now Publishers, Inc.}
}

@article{recht2011hogwild,
  title={{HOGWILD!}: A lock-free approach to parallelizing stochastic gradient descent},
  author={Recht, Benjamin and Re, Christopher and Wright, Stephen and Niu, Feng},
  journal={Advances in Neural Information Processing Systems},
  volume={24},
  year={2011}
}

@article{goyal2017accurate,
  title={Accurate, large minibatch {SGD}: Training imagenet in 1 hour},
  author={Goyal, Priya and Doll{\'a}r, Piotr and Girshick, Ross and Noordhuis, Pieter and Wesolowski, Lukasz and Kyrola, Aapo and Tulloch, Andrew and Jia, Yangqing and He, Kaiming},
  journal={arXiv preprint arXiv:1706.02677},
  year={2017}
}

@article{koloskova2022sharper,
  title={Sharper convergence guarantees for asynchronous {SGD} for distributed and federated learning},
  author={Koloskova, Anastasiia and Stich, Sebastian U and Jaggi, Martin},
  journal={Advances in Neural Information Processing Systems},
  volume={35},
  pages={17202--17215},
  year={2022}
}

@article{cotter2011better,
  title={Better mini-batch algorithms via accelerated gradient methods},
  author={Cotter, Andrew and Shamir, Ohad and Srebro, Nati and Sridharan, Karthik},
  journal={Advances in Neural Information Processing Systems},
  volume={24},
  year={2011}
}

@inproceedings{gower2019sgd,
  title={{SGD}: General analysis and improved rates},
  author={Gower, Robert Mansel and Loizou, Nicolas and Qian, Xun and Sailanbayev, Alibek and Shulgin, Egor and Richt{\'a}rik, Peter},
  booktitle={International Conference on Machine Learning},
  pages={5200--5209},
  year={2019},
  organization={PMLR}
}

@inproceedings{dutta2018slow,
  title={Slow and stale gradients can win the race: Error-runtime trade-offs in distributed {SGD}},
  author={Dutta, Sanghamitra and Joshi, Gauri and Ghosh, Soumyadip and Dube, Parijat and Nagpurkar, Priya},
  booktitle={International Conference on Artificial Intelligence and Statistics},
  pages={803--812},
  year={2018},
  organization={PMLR}
}

@article{chen2016revisiting,
  title={Revisiting distributed synchronous {SGD}},
  author={Chen, Jianmin and Pan, Xinghao and Monga, Rajat and Bengio, Samy and Jozefowicz, Rafal},
  journal={arXiv preprint arXiv:1604.00981},
  year={2016}
}

@inproceedings{nguyen2018sgd,
  title={{SGD} and {Hogwild}! Convergence without the bounded gradients assumption},
  author={Nguyen, Lam and Nguyen, Phuong Ha and Dijk, Marten and Richt{\'a}rik, Peter and Scheinberg, Katya and Tak{\'a}c, Martin},
  booktitle={International Conference on Machine Learning},
  pages={3750--3758},
  year={2018},
  organization={PMLR}
}

@inproceedings{arjevani2020tight,
  title={A tight convergence analysis for stochastic gradient descent with delayed updates},
  author={Arjevani, Yossi and Shamir, Ohad and Srebro, Nathan},
  booktitle={Algorithmic Learning Theory},
  pages={111--132},
  year={2020},
  organization={PMLR}
}

@article{feyzmahdavian2016asynchronous,
  title={An asynchronous mini-batch algorithm for regularized stochastic optimization},
  author={Feyzmahdavian, Hamid Reza and Aytekin, Arda and Johansson, Mikael},
  journal={IEEE Transactions on Automatic Control},
  volume={61},
  number={12},
  pages={3740--3754},
  year={2016},
  publisher={IEEE}
}

@article{mishchenko2022asynchronous,
  title={Asynchronous {SGD} beats minibatch {SGD} under arbitrary delays},
  author={Mishchenko, Konstantin and Bach, Francis and Even, Mathieu and Woodworth, Blake E},
  journal={Advances in Neural Information Processing Systems},
  volume={35},
  pages={420--433},
  year={2022}
}

@article{cohen2021asynchronous,
  title={Asynchronous Stochastic Optimization Robust to Arbitrary Delays},
  author={Cohen, Alon and Daniely, Amit and Drori, Yoel and Koren, Tomer and Schain, Mariano},
  journal={Advances in Neural Information Processing Systems},
  volume={34},
  pages={9024--9035},
  year={2021}
}

@inproceedings{nguyen2022federated,
  title={Federated learning with buffered asynchronous aggregation},
  author={Nguyen, John and Malik, Kshitiz and Zhan, Hongyuan and Yousefpour, Ashkan and Rabbat, Mike and Malek, Mani and Huba, Dzmitry},
  booktitle={International Conference on Artificial Intelligence and Statistics},
  pages={3581--3607},
  year={2022},
  organization={PMLR}
}

@article{tyurin2024optimal,
  title={Optimal time complexities of parallel stochastic optimization methods under a fixed computation model},
  author={Tyurin, Alexander and Richt\'{a}rik, Peter},
  journal={Advances in Neural Information Processing Systems},
  volume={36},
  year={2024}
}

@article{mcmahan2016federated,
  title={Federated learning of deep networks using model averaging},
  author={McMahan, H Brendan and Moore, Eider and Ramage, Daniel and y Arcas, Blaise Ag{\"u}era},
  journal={arXiv preprint arXiv:1602.05629},
  volume={2},
  pages={2},
  year={2016}
}

@inproceedings{islamov2024asgrad,
  title={{AsGrad}: A sharp unified analysis of asynchronous-{SGD} algorithms},
  author={Islamov, Rustem and Safaryan, Mher and Alistarh, Dan},
  booktitle={International Conference on Artificial Intelligence and Statistics},
  pages={649--657},
  year={2024},
  organization={PMLR}
}

@article{tyurin2024shadowheart,
  title={{Shadowheart {SGD}}: Distributed Asynchronous {SGD} with Optimal Time Complexity Under Arbitrary Computation and Communication Heterogeneity},
  author={Tyurin, Alexander and Pozzi, Marta and Ilin, Ivan and Richt{\'a}rik, Peter},
  journal={Advances in Neural Information Processing Systems},
  volume={37},
  year={2024}
}

@article{tyurin2024freya,
  title={Freya {PAGE}: First Optimal Time Complexity for Large-Scale Nonconvex Finite-Sum Optimization with Heterogeneous Asynchronous Computations},
  author={Tyurin, Alexander and Gruntkowska, Kaja and Richt{\'a}rik, Peter},
  journal={Advances in Neural Information Processing Systems},
  volume={37},
  year={2024}
}

@article{tyurin2024optimalgraph,
  title={On the Optimal Time Complexities in Decentralized Stochastic Asynchronous Optimization},
  author={Tyurin, Alexander and Richt{\'a}rik, Peter},
  journal={Advances in Neural Information Processing Systems},
  volume={37},
  year={2024}
}

@InProceedings{yurochkin2019bayesian,
  title = 	 {{B}ayesian Nonparametric Federated Learning of Neural Networks},
  author =       {Yurochkin, Mikhail and Agarwal, Mayank and Ghosh, Soumya and Greenewald, Kristjan and Hoang, Nghia and Khazaeni, Yasaman},
  booktitle = 	 {Proceedings of the 36th International Conference on Machine Learning},
  pages = 	 {7252--7261},
  year = 	 {2019},
  editor = 	 {Chaudhuri, Kamalika and Salakhutdinov, Ruslan},
  volume = 	 {97},
  series = 	 {Proceedings of Machine Learning Research},
  month = 	 {09--15 Jun},
  publisher =    {PMLR},
  url = 	 {https://proceedings.mlr.press/v97/yurochkin19a.html}
}

@article{xiao2017fashion,
  title={Fashion-{MNIST}: a novel image dataset for benchmarking machine learning algorithms},
  author={Xiao, Han and Rasul, Kashif and Vollgraf, Roland},
  journal={arXiv preprint arXiv:1708.07747},
  year={2017}
}

@article{leblond2018improvedAsynchronous,
  author  = {Remi Leblond and Fabian Pedregosa and Simon Lacoste-Julien},
  title   = {Improved Asynchronous Parallel Optimization Analysis for Stochastic Incremental Methods},
  journal = {Journal of Machine Learning Research},
  year    = {2018},
  volume  = {19},
  number  = {81},
  pages   = {1--68},
  url     = {http://jmlr.org/papers/v19/17-650.html}
}

@article{assran2020advances,
  title={Advances in asynchronous parallel and distributed optimization},
  author={Assran, Mahmoud and Aytekin, Arda and Feyzmahdavian, Hamid Reza and Johansson, Mikael and Rabbat, Michael G},
  journal={Proceedings of the IEEE},
  volume={108},
  number={11},
  pages={2013--2031},
  year={2020},
  publisher={IEEE}
}

@inproceedings{zhao2016fast,
  title={Fast asynchronous parallel stochastic gradient descent: A lock-free approach with convergence guarantee},
  author={Zhao, Shen-Yi and Li, Wu-Jun},
  booktitle={Proceedings of the AAAI Conference on Artificial Intelligence},
  volume={30},
  number={1},
  year={2016}
}

@article{reddi2015variance,
  title={On variance reduction in stochastic gradient descent and its asynchronous variants},
  author={J Reddi, Sashank and Hefny, Ahmed and Sra, Suvrit and Poczos, Barnabas and Smola, Alexander J},
  journal={Advances in Neural Information Processing Systems},
  volume={28},
  year={2015}
}

@inproceedings{zhou2018simple,
  title={A simple stochastic variance reduced algorithm with fast convergence rates},
  author={Zhou, Kaiwen and Shang, Fanhua and Cheng, James},
  booktitle={International Conference on Machine Learning},
  pages={5980--5989},
  year={2018},
  organization={PMLR}
}
\bibliographystyle{iclr2026_conference}

\appendix
\newpage

\section{Related Work}
Research on asynchronous stochastic gradient methods dates back to the seminal work of \citet{tsitsiklis1986distributed}, and gained renewed momentum with the introduction of the \algname{HOGWILD!} algorithm \citep{recht2011hogwild}.
\algname{HOGWILD!} is fundamentally an asynchronous coordinate descent method: updates are performed lock-free with inconsistent reads and writes, and its convergence guarantees rely on sparsity assumptions that are rarely realistic in modern large-scale machine learning.
Subsequent refinements of this paradigm include works by \citet{reddi2015variance, zhao2016fast, mania2017perturbed, leblond2018improvedAsynchronous, nguyen2018sgd, zhou2018simple}, but these works remain tied to the coordinate descent setting with inconsistent memory accesses, and thus differ substantially from our focus.

Closer to our setting are works where updates are based on gradients that are applied consistently.
Early contributions, typically under the homogeneous-data assumption (all workers sample from the same distribution), include the work of \citet{agarwal2011distributed}, who studied convex objectives, as well as later extensions to the non-convex case such as the work of \citet{lian2015asynchronous} and \citet{dutta2018slow}, the latter analyzing exponentially distributed computation times.
Other relevant results in this line include those of \citet{feyzmahdavian2016asynchronous, zheng2017asynchronous, arjevani2020tight, feyzmahdavian2023asynchronous}, all of which assume fixed delays.
More recently, delay-adaptive methods have been proposed, aiming to improve performance by down-weighting very stale gradients \citep{cohen2021asynchronous, koloskova2022sharper, mishchenko2022asynchronous}.

Particularly relevant to our work are asynchronous variants of \saga.
\citet{leblond2018improvedAsynchronous} developed a shared-parameter version in the spirit of \algname{HOGWILD!}, while \citet{glasgow2022asynchronous} studied a distributed setting that employs full gradients, in contrast to our stochastic-gradient perspective.

A large body of recent work investigates asynchronous methods in federated learning (FL), where clients hold data from heterogeneous distributions.
Notable contributions include works by \citet{xie2019asynchronous, mishchenko2022asynchronous, koloskova2022sharper, wang2022asyncfeded, wang2022asynchronous, glasgow2022asynchronous, fraboni2023general, zhang2023no, wang2024achieving, islamov2024asgrad, alahyane2025optimizing}.

More broadly, \citet{assran2020advances} provide a comprehensive survey of asynchronous optimization methods.

There is another line of work that began with \citet{tyurin2024optimal}, who established lower bounds for parallel methods and proposed optimal synchronous algorithms together with an asynchronous counterpart.
Several follow-up papers extended this semi-asynchronous framework to other settings \citep{tyurin2024freya,tyurin2024shadowheart,tyurin2024optimalgraph,maranjyan2025mindflayer}.

Finally, beyond asynchronous approaches, several synchronous methods address system heterogeneity by adapting local training to worker speeds.
The canonical method, \algname{FedAvg} \citep{mcmahan2017communication}, performs multiple local steps on each worker.
Variants have adapted the number of local steps to match workers' computation speeds \citep{li2020federated, maranjyan2025gradskip}, effectively balancing task assignments across heterogeneous systems.
More recently, \citet{maranjyan2025ata} proposed adapting the number of local steps dynamically, without requiring prior knowledge of worker speeds.

\section{The Computation-Only Model and the Role of Asynchrony}
\label{sec:communication_limitations}

This section clarifies the scope of our modeling assumptions and the specific phenomenon our results target.
Asynchrony is designed to eliminate \textit{waiting time}: the idle time that arises when workers have heterogeneous computation speeds or experience straggling behavior due to hardware stalls, load imbalance, or even network delays.
A key point is that the goal of asynchrony is \emph{not} to reduce communication cost, but to ensure that slow or delayed workers do not force faster workers to remain idle.

Critically, even in the \emph{simplest} setting---homogeneous data and \emph{no} communication cost---it was unknown until very recently whether an asynchronous SGD method could match the optimal synchronous rate.
In fact, existing asynchronous methods were shown to be \emph{worse} than the optimal synchronous algorithm in this basic regime \citep{tyurin2024optimal}.
The recent work of \citet{maranjyan2025ringmaster} resolved this foundational case for the first time by showing that, under homogeneous data, an asynchronous method can achieve the same optimal time complexity as the synchronous method \rennala \citep{tyurin2024optimal}.
Our work extends this understanding to the significantly more challenging heterogeneous-data setting and shows that asynchrony \emph{can} solve the problem it is designed for---and that it can do so \emph{optimally}.

\paragraph{Asynchrony and stragglers.}
Stragglers may be caused by slow computation, device variability, or even \emph{communication delays}.
From the server's perspective, a worker whose gradient arrives late because it is still communicating appears identical to one that is slow at computing.
Asynchrony ensures that such delays---regardless of source---do not block progress: fast workers continue contributing updates while slow or delayed workers catch up.
This ability to eliminate idle time is precisely the purpose of asynchronous methods.

\paragraph{Asynchrony does \emph{not} reduce communication cost.}
Although asynchrony prevents waiting during communication-induced delays, it does \emph{not} reduce the communication overhead itself.
Reducing communication cost requires \textit{orthogonal techniques} such as gradient compression, sparsification, quantization, or local-update schemes (e.g., \algname{FedAvg}~\citep{mcmahan2016federated}).
For example, \citet{tyurin2024shadowheart} explicitly study communication-aware training and use compression-based methods to reduce communication time---demonstrating that communication efficiency is a \textit{separate algorithmic axis} that must be combined with, rather than replaced by, asynchrony.
Thus, asynchrony resolves the \textit{waiting problem} caused by delays, but not the \textit{communication-cost problem}; addressing the latter requires additional mechanisms.

\paragraph{Why communication is not modeled explicitly here.}
For these reasons, we adopt the standard computation-only model used in essentially all theoretical works on asynchronous SGD  
\citep{mishchenko2022asynchronous, koloskova2022sharper, tyurin2024optimal, tyurin2024optimalgraph,maranjyan2025ringmaster},  
which is also the model under which the lower bounds of \citet{tyurin2024optimal} are derived.
Our claims of \emph{optimal time complexity} therefore refer to this shared and well-established model.
Studying asynchrony under explicit communication cost---where it must interact with compression, local updates, or buffering---requires new lower bounds and a different theoretical framework, and is beyond the scope of this work.

\section{Arbitrarily Changing Computation Times}
\label{sec:arbitrary_time}

In practice, the \emph{fixed computation model} \eqref{eq:fixed_time} is often not satisfied.
The compute power of devices can vary over time due to temporary disconnections, hardware or network delays, fluctuations in processing capacity, or other transient effects \citep{maranjyan2025mindflayer}.

In this section we extend our theory to the more general setting of arbitrarily varying computation times.

\subsection{Universal Computation Model}

To formalize this setting, we adopt the \emph{universal computation model} introduced by \citet{tyurin2024tighttimecomplexitiesparallel}.

For each worker $i \in [n]$, we define a \emph{compute power} function
$$
    p_i : \R_{+} \to \R_{+}~,
$$
assumed nonnegative and continuous almost everywhere (countably many jumps allowed).
For any $T^2 \ge T^1 \ge 0$, the number of stochastic gradients \emph{completed} by worker $i$ on $[T^1, T^2]$ is
$$
    \#\text{gradients in }[T^1, T^2] \;=\; \left\lfloor \int_{T^1}^{T^2} p_i(t)\,dt \right\rfloor.
$$
Here, $p_i(t)$ models the worker's time-varying computational ability: smaller values over an interval yield fewer completed gradients, and larger values yield more.

For instance, if worker $i$ remains idle for the first $T$ seconds and then becomes active, this corresponds to $p_i(t) = 0$ for $t \leq T$ and $p_i(t) > 0$ for $t > T$.
More generally, $p_i(t)$ may follow periodic or irregular patterns, leading to bursts of activity, pauses, or chaotic changes in compute power.
The process $p_i(t)$ may even be random, and all results hold conditional on the realized sample paths of $\{p_i\}$.

The \emph{universal computation model} reduces to the \emph{fixed computation model} \eqref{eq:fixed_time} when $p_i(t) = \nicefrac{1}{\tau_i}$ for all $t \geq 0$ and $i \in [n]$.
In this case,
$$
    \#\text{gradients in }[T^1, T^2] = \left\lfloor \frac{T^2 - T^1}{\tau_i} \right\rfloor,
$$  
meaning that worker $i$ computes one stochastic gradient after $T^1 + \tau_i$ seconds, two gradients after $T^1 + 2\tau_i$ seconds, and so on.

\subsection{Toward an Optimal Method}

In the general setting of arbitrarily varying computation times, \Cref{algo:Ringleader} is not optimal.
To see why, consider the following adversarial timing pattern.

Suppose there are two workers.
During one gradient computation by the slower worker, the faster worker computes $s$ gradients.
Immediately afterwards, they switch roles: the previously fast worker slows down by a factor of $s$, while the previously slow one speeds up by the same factor.
This pattern repeats each time the slower worker finishes a gradient computation.

In this setting, if we run \Cref{algo:Ringleader}, the server waits in each Phase~1 for a single gradient from every worker.
Thus, the slower worker always contributes only one gradient, and the harmonic mean of the batch sizes satisfies 
$$
    1 \;\le\; B^k \;\le\; 2~.
$$ 
From \Cref{theorem:convergence}, the iteration complexity is
$$
    \cO\!\left(
        \frac{nL\Delta}{\varepsilon} \left( 1 + \frac{\sigma^2}{Bn\varepsilon} \right)
    \right).
$$
When $\nicefrac{\sigma^2}{n \varepsilon}$ is much larger than $B$, this dependence can be highly suboptimal.

Instead, suppose the server waits until one full round of the above process completes, collecting $s+1$ gradients from each worker.
Then the harmonic mean satisfies $B^k \ge s+1$, which can be arbitrarily larger than~2.
Since in practice both $s$ and $\nicefrac{\sigma^2}{n\varepsilon}$ can be very large, the naive strategy of waiting for only one gradient per worker (as in \Cref{algo:Ringleader}) cannot be optimal in the arbitrary-time setting.

\subsection{An Optimal Method}

The solution is simple and follows directly from the iteration complexity bound.
From
$$
    \cO\!\left(
        \frac{nL\Delta}{\varepsilon} \left( 1 + \frac{\sigma^2}{Bn\varepsilon} \right)
    \right),
$$
we see that to balance the terms it suffices to ensure
$$
    B \;\ge\; \frac{\sigma^2}{n\varepsilon}~.
$$
Accordingly, we modify the stopping condition in Phase~1 of \Cref{algo:Ringleader}.
Instead of requiring the server to receive at least one gradient from each worker, we require the stronger condition used in \malenia, namely
\begin{equation}\tag{\ref{eq:malenia_condition}}
    \left(\frac{1}{n} \sum_{i=1}^n \frac{1}{b_i} \right)^{-1} 
    \;\;\ge\;\; \max\left\{1, \frac{\sigma^2}{n\varepsilon}\right\}~,
\end{equation}
where $b_i$ is the number of gradients received from worker~$i$.

In the low-noise regime, where $\nicefrac{\sigma^2}{n\varepsilon} \le 1$, the condition reduces to requiring $b_i \ge 1$ for all $i$, so the algorithm coincides with the original \Cref{algo:Ringleader}.
In the high-noise regime, the algorithm collects more gradients in Phase~1, ensuring that $B$ is sufficiently large for optimal convergence.

With this change, Phase~1 of our algorithm matches that of \malenia.
The difference lies in Phase~2: our algorithm continues to use the ongoing gradient computations from all workers to perform $n$ updates, while \malenia discards any unfinished gradients, performs a single update, and then proceeds to the next round.

The following theorem establishes the time complexity of our algorithm under the universal computation model.
\begin{boxedtheorem}
    Under Assumptions~\ref{ass:stochastic_variance_bounded}, \ref{ass:lipschitz_constant}, and \ref{ass:lower_bound}, let the stepsize in \algn be
    $$
        \gamma = \frac{1}{10nL}~.
    $$
    Then, under the \emph{universal computation model}, \algn finds an $\varepsilon$--stationary point within at most $T^{K}$ seconds, where
    $$
        K \;\eqdef\; \left\lceil \frac{160 L \Delta}{\varepsilon} \right\rceil,
    $$
    and $T^{K}$ denotes the $K$-th element of the recursively defined sequence
    \begin{align*}
        T^k \;=\; \min \left\{ T \ge 0 :
            \left( \frac{1}{n}\sum_{i=1}^n 
            \left\lfloor \int_{T_{k-1}}^{T} p_i(t) \, dt \right\rfloor^{-1} \right)^{-1} 
            \;\;\ge\;\; \max\!\left\{1, \frac{\sigma^2}{n\varepsilon}\right\} 
        \right\}~,
    \end{align*}
    for all $k \ge 1$, with initialization $T^0 = 0$.
\end{boxedtheorem}
This result matches the lower bound derived by \citet{tyurin2024tighttimecomplexitiesparallel}, and therefore the proposed method is optimal.  
\begin{proof}
    Under the condition in \eqref{eq:malenia_condition}, each gradient-type step of the algorithm satisfies
    $$
        B^k = \left( \frac{1}{n} \sum_{i=1}^n \frac{1}{b_i^k} \right)^{-1}
        \;\;\ge\;\; \max\left\{1,\frac{\sigma^2}{n\varepsilon}\right\}.
    $$
    In \Cref{theorem:convergence}, instead of using $B$ we can substitute any valid lower bound.
    Here we choose
    $$
        B = \max\left\{1,\frac{\sigma^2}{n\varepsilon}\right\}.
    $$
    With this substitution, the iteration complexity becomes
    $$
        K = \frac{80\,nL\Delta}{\varepsilon}~.
    $$
    To derive the time complexity, consider the time required to perform $n$ iterations.
    Each block of $n$ updates occurs in Phase~2 following the Phase~1 gradient collection.
    Starting from time $T=0$, Phase~1 ends once the accumulated number of gradients satisfies condition~\eqref{eq:malenia_condition}, which occurs at time
    $$
        T_+^1 = \min \left\{T \ge 0 : 
            \left( \frac{1}{n}\sum_{i=1}^n 
            \left\lfloor \int_{0}^{T} p_i(t) \, dt \right\rfloor^{-1} \right)^{-1} 
            \;\;\ge\;\; \max\left\{1, \frac{\sigma^2}{n\varepsilon}\right\} 
        \right\}.
    $$
    After Phase~1, to complete $n$ updates in Phase~2 we must wait for the ongoing computations to finish.
    This requires at most
    $$
        T^1 = \min \left\{T \ge 0 : 
            \left( \frac{1}{n}\sum_{i=1}^n 
            \left\lfloor \int_{T_+^1}^{T} p_i(t) \, dt \right\rfloor^{-1} \right)^{-1} 
            \;\;\ge\;\; 1
        \right\}.
    $$
    Thus, the total time to complete all $K$ iterations is bounded by
    $$
        T^{\left\lceil \nicefrac{2K}{n} \right\rceil}~,
    $$
    where the sequence $\{T^k\}_{k \ge 0}$ is defined recursively as
    $$
        T^k = \min \left\{T \ge 0 : 
            \left( \frac{1}{n}\sum_{i=1}^n 
            \left\lfloor \int_{T_{k-1}}^{T} p_i(t) \, dt \right\rfloor^{-1} \right)^{-1} 
            \;\;\ge\;\; \max\left\{1, \frac{\sigma^2}{n\varepsilon}\right\} 
        \right\}, 
        \qquad T^0 = 0~.
    $$
\end{proof}
\section{Auxiliary Lemmas}
Here we provide proofs of lemmas omitted from the main text, along with auxiliary results that will be used later.
\subsection{Proof of \Cref{lemma:smoothness_relation}}
\label{proof:smoothness}
We begin with a lemma relating the different smoothness constants.
\begin{restate-boxedlemma}{\ref{lemma:smoothness_relation}}[Smoothness Bounds]
    Let $L_f$ denote the smoothness constant of $f$, $L_{f_i}$ the smoothness constant of $f_i$, and $L$ the constant from \Cref{ass:lipschitz_constant}.
    We have
    $$
        L_f \;\le\; L \;\le\; \sqrt{\frac{1}{n}\sum_{i=1}^n L_{f_i}^2} \;\le\; \max_{i \in [n]} L_{f_i} =: L_{\max}~.
    $$
    Moreover, if all $f_i$ are identical, i.e., $f_i = f$ for all $i \in [n]$, then $L = L_f$.
\end{restate-boxedlemma}
Recall from \Cref{ass:lipschitz_constant} that we assumed the following generalized smoothness condition:  
for some constant $L>0$ and for all $x \in \R^d$ and $y_1,\dots,y_n \in \R^d$,  
\begin{align}\label{eq:smoothness}
    \sqnorm{\nabla f(x) - \frac{1}{n}\sum_{i=1}^n \nabla f_i(y_i)}
    \;\le\; \frac{L^2}{n} \sum_{i=1}^n \sqnorm{x - y_i}~.
\end{align}
Recall that a function $\phi$ is called $L_\phi$--smooth if
$$
    \norm{\nabla \phi(x) - \nabla \phi(y)} \le L_\phi \norm{x - y}, \quad \forall x,y \in \R^d~.
$$
Here $L_\phi$ denotes the minimal such constant.
We are ready to prove the lemma.
\begin{proof}
For the first inequality, take $y_1 = \dots = y_n = y$.  
Then \eqref{eq:smoothness} reduces to
$$
    \sqnorm{\nabla f(x) - \nabla f(y)} \le L^2 \sqnorm{x-y},
$$
so $f$ is $L$--smooth.
By definition of $L_f$ as the minimal smoothness constant, this implies $L_f \le L$.

For the second inequality, by the triangle inequality, then by the smoothness of each $f_i$, and finally by Cauchy--Schwarz,
\begin{align*}
    \norm{\nabla f(x) - \frac{1}{n}\sum_{i=1}^n \nabla f_i(y_i)}
        &\le \frac{1}{n}\sum_{i=1}^n \norm{\nabla f_i(x) - \nabla f_i(y_i)}
        \le \frac{1}{n}\sum_{i=1}^n L_{f_i} \norm{x - y_i} \\
        &\le \sqrt{ \frac{1}{n}\sum_{i=1}^n L_{f_i}^2 } \ \sqrt{ \frac{1}{n}\sum_{i=1}^n \sqnorm{x-y_i} }~.
\end{align*}
Squaring both sides shows that \eqref{eq:smoothness} holds with 
$L = \sqrt{ \frac{1}{n}\sum_{i=1}^n L_{f_i}^2 }$~.

Finally, suppose all $f_i$ are identical: $f_i \equiv f$ for all $i$.
Then
\begin{align*}
    \norm{\nabla f(x) - \frac{1}{n}\sum_{i=1}^n \nabla f(y_i)}
    &\le \frac{1}{n}\sum_{i=1}^n \norm{\nabla f(x) - \nabla f(y_i)}
    \le \frac{L_f}{n}\sum_{i=1}^n \norm{x - y_i} \\
    &\le L_f \sqrt{\frac{1}{n}\sum_{i=1}^n \sqnorm{x - y_i}}~,
\end{align*}
where the last step uses Cauchy--Schwarz.
Squaring both sides yields
$$
    \sqnorm{\nabla f(x) - \frac{1}{n}\sum_{i=1}^n \nabla f(y_i)}
    \le \frac{L_f^2}{n}\sum_{i=1}^n \sqnorm{x - y_i},
$$
i.e., \eqref{eq:smoothness} holds with $L \le L_f$.
Combined with $L_f \le L$, we conclude $L = L_f$.
\end{proof}

\subsection{Variance Term}
The following lemma bounds the variance of the gradient estimator in \algn.
\begin{boxedlemma}[Variance Bound]\label{lemma:variance}
    Under \Cref{ass:stochastic_variance_bounded}, the following variance-type inequality holds for the gradient estimator used in \Cref{algo:Ringleader}:
    $$
        \E{\sqnorm{ \bar g^k - \frac{1}{n}\sum_{i=1}^n \nabla f_i\(x^{k-\delta_i^k}\) }}
        \le \frac{\sigma^2}{B^k n}~.
    $$
\end{boxedlemma}
\begin{proof}
    Recall that the gradient estimator is defined as
    $$
        \bar g^k
        = \frac{1}{n} \sum_{i=1}^n \bar g_i^k
        = \frac{1}{n} \sum_{i=1}^n \frac{1}{b_i^k} \sum_{j=1}^{b_i^k} g_i^{k,j}
        = \frac{1}{n} \sum_{i=1}^n \frac{1}{b_i^k} \sum_{j=1}^{b_i^k} \nabla f_i\( x^{k-\delta_i^k}; \xi_i^{k-\delta_i^k, j} \) ~.
    $$
    Let $\mathcal{F}^k$ denote the sigma-field containing all randomness up to the start of the current round, i.e., up to iteration $k-(k \bmod n)$.
    Conditioning on $\mathcal{F}^k$, the evaluation points $x^{k-\delta_i^k}$ are deterministic, and the stochastic gradients $g_i^{k,j}$ are independent across both workers $i$ and samples $j$.

    Using the law of total expectation and the independence of stochastic gradients, we have
    \begin{align*}
        \E{\sqnorm{ \bar g^k - \frac{1}{n}\sum_{i=1}^n \nabla f_i\(x^{k-\delta_i^k}\) }}
        &= \E{\ExpCond{\sqnorm{ \bar g^k - \frac{1}{n}\sum_{i=1}^n \nabla f_i\(x^{k-\delta_i^k}\)}}{\mathcal{F}^k}} \\
        &= \E{\frac{1}{n^2} \sum_{i=1}^{n} \ExpCond{\sqnorm{ \bar g_i^k - \nabla f_i\(x^{k-\delta_i^k}\)}}{\mathcal{F}^k}}.
    \end{align*}
    For each worker $i$, the conditional variance of the minibatch gradient estimator is
    \begin{align*}
        \ExpCond{\sqnorm{ \bar g_i^k - \nabla f_i\(x^{k-\delta_i^k}\)}}{\mathcal{F}^k} 
        &= \ExpCond{\sqnorm{ \frac{1}{b_i^k} \sum_{j=1}^{b_i^k} g_i^{k,j} - \nabla f_i\(x^{k-\delta_i^k}\)}}{\mathcal{F}^k}  \\
        &= \frac{1}{b_i^k} \ExpCond{\sqnorm{ g_i^{k,1} - \nabla f_i\(x^{k-\delta_i^k}\)}}{\mathcal{F}^k}  \le \frac{\sigma^2}{b_i^k}~,
    \end{align*}
    where the last inequality follows from \Cref{ass:stochastic_variance_bounded}.
    
    Combining these results, we get
    \begin{align*}
        \E{\sqnorm{ \bar g^k - \frac{1}{n}\sum_{i=1}^n \nabla f_i\(x^{k-\delta_i^k}\) }}
        &\le \frac{1}{n^2} \sum_{i=1}^{n} \frac{\sigma^2}{b_i^k}
            = \frac{\sigma^2}{n} \frac{1}{n} \sum_{i=1}^n \frac{1}{b_i^k} 
            = \frac{\sigma^2}{B^k n}~,
    \end{align*}
    where the last equality uses the definition of the harmonic mean $B^k = \(\frac{1}{n} \sum_{i=1}^n \frac{1}{b_i^k}\)^{-1}$.
\end{proof}
\subsection{Proof of \Cref{lemma:descent}}
\label{proof:descent}
We now prove the descent lemma.
\begin{restate-boxedlemma}{\ref{lemma:descent}}[Descent Lemma]
    Under Assumptions~\ref{ass:stochastic_variance_bounded} and~\ref{ass:lipschitz_constant}, if the stepsize in \Cref{algo:Ringleader} satisfies $\gamma \le \nicefrac{1}{4L}$, then the following inequality holds
    \begin{align*}
        \E{f\(x^{k+1}\)}
        &\le \E{f\(x^{k}\)}
            - \frac{\gamma}{2} \E{ \sqnorm{\nabla f\(x^{k}\)} }
            - \frac{\gamma}{4} \E{\sqnorm{\frac{1}{n} \sum_{i=1}^n \nabla f_i\(x^{k-\delta_i^k}\)}} \\
            &\quad + \frac{\gamma L^2}{2n} \sum_{i=1}^n \E{\sqnorm{x^{k} - x^{k-\delta_{i}^k}}}
                    + \frac{3\gamma^2 L \sigma^2}{2B} \\
            &\quad + \gamma^2 L \sum_{\ell = k-(k \bmod n)}^{k-1} \E{\sqnorm{\frac{1}{n}\sum_{i=1}^n \nabla f_i\(x^{\ell-\delta_i^\ell}\)}}.
    \end{align*}
\end{restate-boxedlemma}
\begin{proof}
    Some proof techniques are adapted from the works of \citet{maranjyan2025ringmaster} and \citet{wang2025incremental}.

    From \Cref{ass:lipschitz_constant} and \Cref{lemma:smoothness_relation}, we know that $f$ is $L$--smooth. 
    Therefore, the following standard inequality holds \citep{nesterov2018lectures}
    \begin{equation}\label{eq:l-smoothness}
        \E{f(x^{k+1})}
        \le \E{ f(x^{k}) 
            - \gamma \inp{\nabla f(x^{k})}{\bar g^k}
            + \frac{L\gamma^2}{2} \sqnorm{\bar g^k} }.
    \end{equation}
    Recall that the gradient estimator is defined as
    $$
        \bar g^k
        = \frac{1}{n} \sum_{i=1}^n \bar g_i^k
        = \frac{1}{n} \sum_{i=1}^n \frac{1}{b_i^k} \sum_{j=1}^{b_i^k} g_i^{k,j} ~.
    $$
    Let $\mathcal{F}^k$ denote the sigma-field containing all randomness up to the start of the current Phase~2, i.e., up to iteration $k - (k \bmod n)$.
    A key observation is that all gradients in the current gradient table were computed and received during the current round.
    Since these gradients were computed at points from previous iterations within the current round, we have $k - \delta_i^k \le k - (k \bmod n)$ for all $i\in[n]$.
    Conditioning on $\mathcal{F}^k$, the points $x^{k-\delta_i^k}$ are deterministic.
    Therefore, we can compute the conditional expectation of the gradient estimator:
    $$
        \ExpCond{\bar g^k}{\mathcal{F}^k}
        = \frac{1}{n} \sum_{i=1}^n \frac{1}{b_i^k} \sum_{j=1}^{b_i^k} \ExpCond{g_i^{k,j}}{\mathcal{F}^k}
        = \frac{1}{n} \sum_{i=1}^n \nabla f_i\(x^{k-\delta_i^k}\).
    $$
    The last equality follows from the unbiasedness of the stochastic gradient estimator (\Cref{ass:stochastic_variance_bounded}).
    
    Using this conditional expectation and the law of total expectation, we can now simplify the inner product term in \eqref{eq:l-smoothness}:
    \begin{align*}
        \E{\inp{\nabla f\(x^{k}\)}{\bar g^k}}
        &= 
        \E{\inp{\nabla f\(x^{k}\)}{\frac{1}{n} \sum_{i=1}^n \nabla f_i\(x^{k-\delta_i^k}\)}} \\
            &\quad + \E{\inp{\nabla f\(x^{k}\)}
            {\bar g^k - \frac{1}{n} \sum_{i=1}^n  \nabla f_i\(x^{k-\delta_i^k}\) }} \\
        &= 
        \E{\inp{\nabla f\(x^{k}\)}{\frac{1}{n} \sum_{i=1}^n \nabla f_i\(x^{k-\delta_i^k}\)}} \\
            &\quad + \E{\inp{\nabla f\(x^{k}\) - \nabla f\(x^{k-(k \bmod n)}\)}
            {\bar g^k - \frac{1}{n} \sum_{i=1}^n  \nabla f_i\(x^{k-\delta_i^k}\) }} \\
                &\qquad + \E{\inp{\nabla f\(x^{k-(k \bmod n)}\)}
                {\bar g^k - \frac{1}{n} \sum_{i=1}^n \nabla f_i\(x^{k-\delta_i^k}\) }} \\
        &= 
        \underbrace{ \E{\inp{\nabla f\(x^{k}\)}{\frac{1}{n} \sum_{i=1}^n \nabla f_i\(x^{k-\delta_i^k}\)}} }_{T_1} \\
            &\quad + \underbrace{ \E{\inp{\nabla f\(x^{k}\) - \nabla f\(x^{k-(k \bmod n)}\)}
            {\bar g^k - \frac{1}{n} \sum_{i=1}^n  \nabla f_i\(x^{k-\delta_i^k}\) }} }_{T_2} .
    \end{align*}
    Next, using \Cref{ass:lipschitz_constant}, we have
    \begin{align*}
        2T_1
        &= \E{2 \inp{\nabla f\(x^{k}\)}{ \frac{1}{n} \sum_{i=1}^n \nabla f_i\(x^{k-\delta_i^k}\) }} \\
        &= \E{ \sqnorm{\nabla f\(x^{k}\)} 
            + \sqnorm{ \frac{1}{n} \sum_{i=1}^n \nabla f_i\(x^{k-\delta_i^k}\) } } 
                - \E{\sqnorm{\nabla f\(x^{k}\) - \frac{1}{n} \sum_{i=1}^n \nabla f_i\(x^{k-\delta_i^k}\) }} \\
        &\ge \E{ \sqnorm{\nabla f\(x^{k}\)} }
            + \E{ \sqnorm{ \frac{1}{n} \sum_{i=1}^n \nabla f_i\(x^{k-\delta_i^k}\) } }
                - \frac{L^2 }{n} \sum_{i=1}^n \E{\sqnorm{x^{k} - x^{k-\delta_i^k}}}.
    \end{align*}
    Next, we analyze $T_2$
    \begin{align*}
        T_2
        &= \E{\inp{\nabla f\(x^{k}\) - \nabla f\(x^{k-(k \bmod n)}\)}
        {\bar g^k - \frac{1}{n} \sum_{i=1}^n  \nabla f_i\(x^{k-\delta_i^k}\) }} \\
        &\ge - \E{\norm{\nabla f\(x^{k}\) - \nabla f\(x^{k-(k \bmod n)}\)} 
            \norm{\bar g^k - \frac{1}{n} \sum_{i=1}^n  \nabla f_i\(x^{k-\delta_i^k}\)}} \\
        &\ge - L \E{\norm{x^{k} - x^{k-(k \bmod n)} }
            \norm{\bar g^k - \frac{1}{n} \sum_{i=1}^n  \nabla f_i\(x^{k-\delta_i^k}\)}} \\
        &= - L \E{\norm{\gamma \sum_{\ell = k-(k \bmod n)}^{k-1} \bar g^\ell }
            \norm{\bar g^k - \frac{1}{n} \sum_{i=1}^n  \nabla f_i\(x^{k-\delta_i^k}\)}} \\
        &\ge - L\gamma \sum_{\ell = k-(k \bmod n)}^{k-1} \E{\norm{ \bar g^\ell }
            \norm{\bar g^k - \frac{1}{n} \sum_{i=1}^n  \nabla f_i\(x^{k-\delta_i^k}\)}} \\
        &\ge - L\gamma \sum_{\ell = k-(k \bmod n)}^{k-1} \frac{1}{2}\(\E{\sqnorm{ \bar g^\ell }} 
            + \E{\sqnorm{\bar g^k - \frac{1}{n} \sum_{i=1}^n  \nabla f_i\(x^{k-\delta_i^k}\)}} \) \\
        &\ge - \frac{L\gamma}{2} \sum_{\ell = k-(k \bmod n)}^{k-1} \E{\sqnorm{ \bar g^\ell } }
            - (k \bmod n) \frac{L\gamma \sigma^2}{2B^k n} \\
        &\ge - \frac{L\gamma}{2} \sum_{\ell = k-(k \bmod n)}^{k-1} \E{\sqnorm{ \bar g^\ell } }
            - \frac{L\gamma \sigma^2}{2B^k}~.
    \end{align*}
    The inequalities follow from the Cauchy-Schwarz inequality, $L$--smoothness of $f$, the triangle inequality, Young's inequality, \Cref{lemma:variance}, and finally $(k \bmod n) \le n-1 < n$.

    It remains to bound the term $\E{\sqnorm{ \bar g^k } }$.
    Using Young's inequality, we have
    \begin{align*}
        \E{\sqnorm{\bar g^k}} 
        &\le 2\E{\sqnorm{\frac{1}{n}\sum_{i=1}^n \nabla f_i\(x^{k-\delta_i^k}\)}} 
            + 2\E{\sqnorm{ \bar g^k - \frac{1}{n}\sum_{i=1}^n \nabla f_i\(x^{k-\delta_i^k}\)}} \\
        &\le 2\E{\sqnorm{\frac{1}{n}\sum_{i=1}^n \nabla f_i\(x^{k-\delta_i^k}\)}} 
            + \frac{2 \sigma^2}{B^k n}~,
    \end{align*}
    where in the last step we used \Cref{lemma:variance}.

    Now, by combining all terms in \eqref{eq:l-smoothness}, we obtain
    \begin{align*}
        \E{f\(x^{k+1}\)}
        &\le \E{f\(x^{k}\)}
            - \frac{\gamma}{2} \E{ \sqnorm{\nabla f\(x^{k}\)} }
            - \frac{\gamma}{2} \E{\sqnorm{\frac{1}{n} \sum_{i=1}^n \nabla f_i\(x^{k-\delta_i^k}\)}} \\
            &\quad + \frac{\gamma L^2}{2n} \sum_{i=1}^n \E{\sqnorm{x^{k} - x^{k-\delta_{i}^k}}} \\
            &\quad + \frac{\gamma^2 L}{2} \sum_{\ell = k-(k \bmod n)}^{k-1} \E{\sqnorm{ \bar g^\ell } }
            + \frac{\gamma^2 L \sigma^2}{2 B^k} + \frac{\gamma^2L}{2} \E{\sqnorm{\bar g^k}} \\
        &\le \E{f\(x^{k}\)}
            - \frac{\gamma}{2} \E{ \sqnorm{\nabla f\(x^{k}\)} }
            - \frac{\gamma}{2} \E{\sqnorm{\frac{1}{n} \sum_{i=1}^n \nabla f_i\(x^{k-\delta_i^k}\)}} \\
            &\quad + \frac{\gamma L^2}{2n} \sum_{i=1}^n \E{\sqnorm{x^{k} - x^{k-\delta_{i}^k}}} \\
            &\quad + \gamma^2 L \sum_{\ell = k-(k \bmod n)}^{k-1} \E{\sqnorm{\frac{1}{n}\sum_{i=1}^n \nabla f_i\(x^{\ell-\delta_i^\ell}\)}} \\
            &\quad + \gamma^2L \E{\sqnorm{\frac{1}{n}\sum_{i=1}^n \nabla f_i\(x^{k-\delta_i^k}\)}} \\
            &\quad + \gamma^2 L \sum_{\ell = k-(k \bmod n)}^{k-1}\frac{\sigma^2}{B^\ell n} 
                + \frac{\gamma^2 L \sigma^2 }{2 B^k} 
                + \frac{\gamma^2L \sigma^2}{B^k n} ~.
    \end{align*}
    This completes the proof under the stepsize condition $\gamma \le \nicefrac{1}{4L}$ and $B \eqdef \inf_{k \ge 0} B^k$.
\end{proof}
\subsection{Proof of \Cref{lemma:residual}}
\label{proof:residual}
The following lemma provides an upper bound on the residual error due to delays.
\begin{restate-boxedlemma}{\ref{lemma:residual}}[Residual Estimation]
    Under \Cref{ass:stochastic_variance_bounded}, the iterates of \algn (\Cref{algo:Ringleader}) with stepsize $\gamma \le \nicefrac{1}{4nL}$ satisfy the following bound
    \begin{equation*}
    \frac{1}{K} \sum_{k=0}^{K-1} \frac{1}{n} \sum_{i=1}^n \E{\sqnorm{x^{k} - x^{k-\delta_i^k}}}
        \le \frac{2\gamma n}{LK}\sum_{k=0}^{K-1} \E{\sqnorm{ \frac{1}{n}\sum_{j=1}^n \nabla f_j\(x^{k-\delta_j^k}\) }}
            + \frac{2\gamma \sigma^2}{LB}~.
    \end{equation*}
\end{restate-boxedlemma}
\begin{proof}
    By Young's inequality, we have
    \begin{align*}
    \E{\sqnorm{x^{k} - x^{k-\delta_i^k}}}
        &= \E{\sqnorm{\gamma \sum_{\ell = k-\delta_i^k}^{k-1} \bar g^\ell}} \\
        &\le 2\gamma^2 \E{\sqnorm{\sum_{\ell = k-\delta_i^k}^{k-1} \frac{1}{n}\sum_{j=1}^n \nabla f_j\(x^{\ell-\delta_j^\ell}\) }} \\
        &\quad + 2\gamma^2 \E{\sqnorm{\sum_{\ell = k-\delta_i^k}^{k-1} \( \bar g^\ell - \frac{1}{n}\sum_{j=1}^n \nabla f_j\(x^{\ell-\delta_j^\ell}\) \)}} \\
        &\le 2\gamma^2 \underbrace{\delta_i^k \sum_{\ell = k-\delta_i^k}^{k-1} \E{\sqnorm{ \frac{1}{n}\sum_{j=1}^n \nabla f_j\(x^{\ell-\delta_j^\ell}\) }}}_{T_{ik}} + 2\gamma^2 (\delta_i^k)^2 \frac{\sigma^2}{Bn}~.
    \end{align*}
    In the last inequality, we used Jensen's inequality and \Cref{lemma:variance}.
    
    Next, we estimate the sum of $T_{ik}$
    \begin{align*}
    \sum_{k=0}^{K-1} \frac{1}{n} \sum_{i=1}^n T_{ik}
        &= \sum_{k=0}^{K-1} \frac{1}{n} \sum_{i=1}^n \delta_i^k \sum_{\ell = k-\delta_i^k}^{k-1} \E{\sqnorm{ \frac{1}{n}\sum_{j=1}^n \nabla f_j\(x^{\ell-\delta_j^\ell}\) }} \\
        &= \frac{1}{n} \sum_{i=1}^n \sum_{k=0}^{K-1} \delta_i^k \sum_{\ell = k-\delta_i^k}^{k-1} \E{\sqnorm{ \frac{1}{n}\sum_{j=1}^n \nabla f_j\(x^{\ell-\delta_j^\ell}\) }} \\
        &\le 2 \sum_{i=1}^n \sum_{k=0}^{K-1} \sum_{\ell = k-\delta_i^k}^{k-1} \E{\sqnorm{ \frac{1}{n}\sum_{j=1}^n \nabla f_j\(x^{\ell-\delta_j^\ell}\) }} \\
        &\le 2 \sum_{i=1}^n \delta_i^{\max} \sum_{k=0}^{K-1} \E{\sqnorm{ \frac{1}{n}\sum_{j=1}^n \nabla f_j\(x^{k-\delta_j^k}\) }} \\
        &\le 4n^2 \sum_{k=0}^{K-1} \E{\sqnorm{ \frac{1}{n}\sum_{j=1}^n \nabla f_j\(x^{k-\delta_j^k}\) }} .
    \end{align*}
    In the first and last inequality, we used the bound $\delta_i^{\max} \le 2n$ from \Cref{lemma:delay}.
    Finally, applying the stepsize condition $\gamma \le \nicefrac{1}{4nL}$ yields the result.
\end{proof}

\section{Time Complexity of \iasgdtitle}
\label{sec:iasgd_time_complexity}

The iteration complexity of \iasgd \citep{wang2025incremental} is
$$
    K = \cO\!\left(\frac{\delta^{\max}L\Delta}{\varepsilon} \left(1 + \frac{\sigma^2}{n\varepsilon}\right)\right).
$$
We now analyze the corresponding wall-clock time under the \textit{fixed computation model} \eqref{eq:fixed_time}.
Since the algorithm performs an update whenever a single worker finishes a computation, we seek the minimal time $T$ such that
$$
    \sum_{i=1}^n\left\lfloor\frac{T}{\tau_i}\right\rfloor \;\ge\; K~.
$$
Observe that
$$
    \sum_{i=1}^n \frac{T}{\tau_i}
        \;\ge\; \sum_{i=1}^n\left\lfloor\frac{T}{\tau_i}\right\rfloor.
$$
Hence, if we define $T^\prime$ by
$$
    \sum_{i=1}^n \frac{T^\prime}{\tau_i} \;=\; K~,
$$
then 
$$
    T^\prime = \left(\sum_{i=1}^n \frac{1}{\tau_i}\right)^{-1} K~.
$$
It follows that the minimal time $T$ is necessarily larger than $T^\prime$.

It remains to bound $\delta^{\max}$.
At initialization, all workers start computing their first gradients simultaneously.
By the time the slowest worker completes its first gradient (at time $\tau_n$), the other workers may each have completed multiple gradients.
In particular,
$$
    \delta^{\max} \;\ge\; \sum_{i=1}^n \left\lfloor \frac{\tau_n}{\tau_i}\right\rfloor.
$$
Combining this with the iteration complexity bound, we obtain that the total runtime satisfies
$$
    T \;\ge\; c \times \frac{\tau_n L\Delta}{\varepsilon} \left(1 + \frac{\sigma^2}{n\varepsilon}\right),
$$
for some universal constant $c > 0$.

Note that the expression above should not be viewed as an exact upper bound on the runtime.
It is better understood as a simplified estimate of~$T$, which is sufficient for our purposes
and provides a cleaner basis for comparison.

\section{Improved Malenia SGD}\label{sec:malenia_param_free}
\malenia has the following iteration complexity \citep{tyurin2024optimal}
$$
    K \;\ge\; \frac{12\Delta L_f}{\varepsilon} \;+\; \frac{12 \Delta L_f \sigma^2}{\varepsilon^2 n S}~,
$$
where $S$ is a lower bound on the harmonic mean of the batch sizes, i.e.,
$$
    \left(\frac{1}{n}\sum_{i=1}^n \frac{1}{b_i^k}\right)^{-1} \;\ge\; S~,
$$
for all iterations $k$.  
In the original \malenia analysis \citep{tyurin2024optimal}, this bound follows from the condition in \eqref{eq:malenia_condition}, which fixes the same value of $S$ across all iterations.

In the fixed-time regime \eqref{eq:fixed_time}, however, this condition is no longer necessary.  
By adopting the same strategy as \algn (\Cref{algo:Ringleader})---namely, waiting for at least one gradient from each worker---we effectively replace $S$ with $B$ in the rate.  
This yields the following time complexity
$$
    \tau_n K
    \;=\; \frac{12\tau_n\Delta L_f}{\varepsilon}\left(1 + \frac{\sigma^2}{\varepsilon n B}\right).
$$
Substituting the expression for $B$ from \Cref{lemma:time_for_n_iter} and proceeding as in the proof of \Cref{thm:time_complexity}, we obtain the same overall time complexity as before---this time \emph{without} requiring condition \eqref{eq:malenia_condition}, which depends on knowing $\sigma$ and fixing $\varepsilon$ in advance.

Finally, note that this improvement is only valid in the fixed-time regime.  
In the setting with arbitrarily varying computation times, the same optimization cannot be applied, for the same reasons discussed for \algn in \Cref{sec:arbitrary_time}.
\end{document}